\documentclass[12pt,a4paper]{article}

\usepackage{latexsym,epsfig}
\usepackage{graphicx} 
\usepackage[usenames,dvipsnames]{color} 
\usepackage{url} 
\usepackage{soul} 
\usepackage{array} 
\usepackage{amssymb}
\usepackage{amsfonts}
\usepackage{amsmath}
\usepackage{hyperref}
\usepackage{verbatim}
\usepackage{anysize}
\marginsize{1cm}{1cm}{1cm}{1cm}

\newtheorem{lemma}{Lemma}[section]
\newtheorem{proposition}{Proposition}[section]
\newtheorem{theorem}{Theorem}[section]
\newtheorem{corollary}{Corollary}[section]

\newtheorem{conjecture}{Conjecture}[section]
\newtheorem{observation}{Observation}[section]

\newcommand{\EndProof}{\hspace{\stretch{1}} $\Box$}

\newcommand{\pr}{\noindent{\bf Proof.}\ }
\newcommand{\G}{\Gamma}
\newcommand{\al}{\alpha}
\newcommand{\D}{\Delta}
\newcommand{\C}{\mathcal{C}}
\newcommand{\rep}{\mathrm{rep}}
\newcommand{\M}{\mathrm{M}}
\newcommand{\N}{\mathrm{N}}
\newcommand{\g}{\mathrm{g}}
\newcommand{\m}{\mathrm{m}}

\title{On bipartite graphs of defect at most 4}
\author{Ramiro Feria-Pur\'on\footnote{\href{mailto:ramiro.feria.puron@gmail.com}{ramiro.feria.puron@gmail.com}} \vspace{-1.5mm}\\
\emph{\small  School of Electrical Engineering and Computer Science} \vspace{-1.7mm}\\{\it \small The University of Newcastle, Australia}\and
Guillermo
Pineda-Villavicencio\footnote{Corresponding author: \href{mailto:work@guillermo.com.au}{work@guillermo.com.au}} \vspace{-1.5mm}
\\\emph{\small Centre for Informatics and Applied Optimisation} \vspace{-1.7mm}\\
\emph{\small University of Ballarat, Australia}}
\begin{document}
\maketitle

\begin{abstract}
\noindent We consider the bipartite version of the {\it degree/diameter problem}, namely, given natural numbers $\Delta\ge2$ and $D\ge2$, find the maximum number $\mathrm{N}^b(\Delta,D)$ of vertices  in a bipartite graph of maximum degree $\Delta$ and diameter $D$. In this context, the Moore bipartite bound $\mathrm{M}^b(\Delta,D)$ represents an upper bound for $\mathrm{N}^b(\Delta,D)$.

\noindent Bipartite graphs of maximum degree $\Delta$, diameter $D$ and order $\mathrm{M}^b(\Delta,D)$ -- called {\it Moore bipartite graphs} -- turned out to be very rare. Therefore, it is very interesting to investigate bipartite graphs of maximum degree $\Delta\ge2$, diameter $D\ge2$  and order $\mathrm{M}^b(\Delta,D)-\epsilon$ with small $\epsilon> 0$; that is, bipartite $(\Delta,D,-\epsilon)$-graphs. The parameter $\epsilon$ is called the {\it defect}.

\noindent This paper considers bipartite graphs of defect at most $4$, and presents all the known such graphs. Bipartite graphs of defect $2$ have been studied in the past; if $\Delta\ge3$ and $D\ge3$, they may only exist for $D=3$. However, when $\epsilon>2 $ bipartite $(\Delta,D,-\epsilon)$-graphs represent a wide unexplored area.

\noindent The main results of the paper include several necessary conditions for the existence of bipartite $(\Delta,d,-4)$-graphs; the complete catalogue of bipartite $(3,D,-\epsilon)$-graphs with $D\ge 2$ and $0\le \epsilon \le4$;  the complete catalogue of bipartite $(\Delta,D,-\epsilon)$-graphs with $\Delta\ge 2$, $5\le D\le 187$ ($D\ne6$) and  $0\le \epsilon \le4$; and a non-existence proof of all bipartite $(\Delta,D,-4)$-graphs with $\Delta\ge 3$ and odd $D\ge 5$.

\noindent Finally, we conjecture that there are no bipartite graphs of defect $4$ for $\Delta\ge3$ and $D\ge5$, and comment on some implications of our results for the upper bounds of $\mathrm{N}^b(\Delta,D)$.
\end{abstract}




\section{Introduction}

Due to the diverse features and applications of interconnection networks, it is possible to find many interpretations about network ``optimality'' in the literature. Here we are concerned with the following; see {\cite[pp. 18]{DYN03}}, {\cite[pp. 168]{Hey96}}, and {\cite[pp. 91 ]{Xu01}}.

\begin{quote}
{\it An optimal network contains the maximum possible
number of nodes, given a limit on the number of connections
attached to a node and a limit on the distance between any two nodes of the network.}
\end{quote}

This interpretation has attracted network designers and the research community in general due to its implications in the design of large interconnection networks. In graph-theoretical terms, this interpretation leads to the {\it degree/diameter
problem} (the problem of finding the largest possible number of vertices in a graph with given  maximum degree and diameter). If the graphs in question are subject to further restrictions such as being bipartite, planarity and/or transitivity, we can state the degree/diameter problem for the classes of graphs under consideration.

In this paper we will consider only bipartite graphs, and in this case, the degree/diameter problem can be stated as follows.
\begin{itemize}
\item[]{\it Degree/diameter problem for bipartite graphs}: Given natural numbers $\Delta\ge2$
and $D\ge2$, find the largest possible number $\N^b(\D,D)$ of vertices in a bipartite graph of maximum degree $\Delta$ and diameter $D$.
\end{itemize}
Note that $\N^b(\D,D)$ is well defined for $\D\ge2$ and $D\ge2$. An upper bound for $\N^b(\D,D)$ is given by the {\it Moore bipartite bound} $\M^b(\D,D)$, defined below:
\[
\M^b(\Delta,D)=2\left(1+(\Delta-1)+\dots+(\Delta-1)^{D-1}\right).
\]
Bipartite graphs of degree $\D$, diameter $D$ and  order $\M^b(\Delta,D)$ are called {\it Moore bipartite graphs}. Moore bipartite graphs are rare; for $\D=2$ they are the cycles of length $2D$, while for $\D\ge3$ Moore bipartite graphs exist only for diameters 2, 3, 4 and 6; see \cite{FH64}. Therefore, we are interested in studying the existence or otherwise of bipartite graphs of given maximum degree $\Delta$, diameter $D$ and order $\M^b(\Delta,D)-\epsilon$ for $\epsilon>0$; that is, bipartite $(\Delta,D,-\epsilon)$-graphs, where the parameter $\epsilon$ is called the \emph{defect}. For notational convenience, we consider Moore bipartite graphs as having defect $\epsilon=0$.

Only a few values of $\N^b(\D,D)$ are known at present. With the exception of $\N^b(3,5)=\M^b(3,5)-6$, settled in \cite{Jor93}, the other known values of $\N^b(\D,D)$ are those for which there is a Moore bipartite graph. The paper \cite{PV} combined with \cite{DJMP2,DJMP3} almost settled the case of bipartite graphs of defect 2; if $\D\ge 3$ and $D\ge3$, then such graphs may only exist for  $D=3$ and certain values of $\D$. Bipartite $(\Delta,D,-\epsilon)$-graphs with $\epsilon>2$ have been rarely considered in the literature so far.

In this paper we consider bipartite $(\D,D,-4)$-graphs with $\D\ge 2$ and $D\ge 3$. By using combinatorial approaches we obtain several important results about bipartite graphs of defect 4, including several necessary conditions for the existence of bipartite $(\Delta,d,-4)$-graphs; the complete catalogue of bipartite $(3,D,-\epsilon)$-graphs with $D\ge 2$ and $0\le \epsilon \le4$;  the complete catalogue of bipartite $(\D,D,-\epsilon)$-graphs with $\D\ge 2$, $5\le D\le 187$ ($D\ne6$) and $0\le \epsilon \le4$; and a non-existence proof of all bipartite $(\D,D,-4)$-graphs with $\D\ge 3$ and odd $D\ge 5$. Finally, we conjecture that there are no bipartite graphs of defect $4$ for $\D\ge3$ and $D\ge5$.

The main results in this paper do not apply to bipartite $(\D,D,-4)$-graphs with $\D\ge 4$ and $D=3,4$. Some of our assertions, however, do offer a partial characterisation of all bipartite $(\D,D,-4)$-graphs with $\D\ge 3$ and $D\ge3$. At the time of writing the paper we do not foresee a conclusive way to take on the diameters $3$ and $4$. To deal with such graphs it would be necessary to either  find different ideas or  complement some of the ones presented here. Section \ref{sub:Diameters34} contains further comments on such diameters.


\section{Notation and Terminology}
\label{sec:Notation}
The terminology and notation used in this paper is standard and consistent with that used in \cite{Die05}, so only those concepts that can vary from texts to texts will be defined.

All graphs considered are simple. The vertex set of a graph $\G$ is denoted by $V(\G)$, and its edge set by $E(\G)$. The \emph{difference} between the graphs $\G$ and $\G'$, denoted by $\G-\G'$, is the graph with vertex set $V({\G})-V({\G'})$ and edge set formed by all the edges with both endvertices in $V({\G})-V({\G'})$.

The set of neighbours of a vertex $x$ in $\G$ is denoted by $N(x)$. For an edge $e=\{x,y\}$ we write $e=xy$, or alternatively $x\sim y$. The set of edges in a graph $\G$ joining a vertex $x$ in $X\subseteq V(\G)$ to a vertex $y$ in $Y\subseteq V(\G)$ is denoted by $E(X,Y)$; for simplicity, we write $E(x,Y)$ rather than  $E(\{x\},Y)$.

A path of length $k$ is called a {\it $k$-path}, and cycle of length $k$ is called a {\it $k$-cycle}. A path from a vertex $x$ to a vertex $y$ is denoted by $x-y$. Whenever we refer to paths we mean shortest paths. We will use the following notation for subpaths of a path $P=x_0x_1\ldots x_k$: $x_iPx_j = x_i\ldots x_j$, where $0\leq i\leq j\leq k$. The distance between a vertex $x$ and a vertex $y$ is denoted by $d(x,y)$.

The union of three independent paths of length $D$ with common endvertices is denoted by $\Theta_D$. In a graph $\G$ a vertex of degree at least 3 is called a \emph{branch vertex} of $\G$.


\section{Known bipartite $(\D,D,-\epsilon)$-graphs with $\D\ge2$, $D\ge2$ and $0\le\epsilon\le4$}
\label{sec:KnownGraphs}
For $\D=2$ the Moore bipartite graphs are the cycles on $2D$ vertices, while for $D=2$ and each $\D\ge 3$ they are the complete bipartite graphs of degree $\D$ and order $2\D$. For $D=3,4,6$ Moore bipartite graphs of degree $\D$ have  been constructed only when $\D-1$ is a prime power \cite{Ben66}. Furthermore, Singleton \cite{Si66} proved that the existence of a Moore bipartite graph of diameter $3$ is equivalent to the existence of a projective plane of order $\D-1$. The question of whether Moore bipartite graphs of diameter $3$, $4$ or $6$ exist for other values of $\D$ remains open, and represents one of the most famous problems in combinatorics. For other values of $D\ge2$ and $\D\ge3$ there are no Moore bipartite graphs (see \cite{Si66, FH64}).

When $\D=2$ or $D=2$ bipartite $(\D,D,-\epsilon)$-graphs with $\epsilon\ge1$ can be obtained by simple observation. For a given $D\ge 2$ there is only one bipartite $(2,D,-\epsilon)$-graph with $\epsilon\ge1$: the path of length $D$, which has defect $\epsilon=D-1$. For a given $\D\ge 2$ there are exactly $\D-1$ bipartite $(\D,2,-\epsilon)$-graphs with $\epsilon\ge1$; they are the complete bipartite graphs with partite sets of size $\D$ and $\D-\epsilon$, where $1 \le \epsilon \le\D-1$. Therefore, from now on we assume $\D\ge3$ and $D\ge 3$.

We continue with some conditions for the regularity of
bipartite $(\Delta,D,-\epsilon)$-graphs, which were obtained in \cite{DJMP2}.

\begin{proposition}[\cite{DJMP2}]
\label{prop:regularity1} For $\epsilon<
1+(\Delta-1)+(\Delta-1)^2+\ldots+(\Delta-1)^{D-2}$, $
\Delta\ge 3$ and $D\ge 3$, a
bipartite $(\Delta,D,-\epsilon)$-graph is regular. 
\end{proposition}

\begin{proposition}[\cite{DJMP2}]
\label{prop:regularity2} For $\epsilon<
2\left((\Delta-1)+(\Delta-1)^3+\ldots+(\Delta-1)^{D-2}\right)$,
$\Delta\ge 3$ and odd $D\ge 3$, a
bipartite $(\Delta,D,-\epsilon)$-graph is regular. 

\end{proposition}

By Propositions \ref{prop:regularity1} and \ref{prop:regularity2}, bipartite $(\D,D,-\epsilon)$-graphs with $\D\ge 3$,  $D\ge 3$ and $\epsilon\le3$ must be regular, implying the non-existence of  such graphs for $\D\ge 3$,  $D\ge 3$ and $\epsilon=1,3$. In the same way, bipartite $(\D,D,-4)$-graphs with $\D\ge 3$ and $D\ge 4$ and bipartite $(\D,3,-4)$-graphs with $\D\ge 4$ must be regular.

For $\D\ge3$ and $D\ge3$, the only known bipartite $(\D,D,-2)$-graphs are depicted in Fig.~\ref{fig:Bipartite(d,3,-2)}. Recall that such graphs do not exist when $\D\ge3$ and $D\ge4$; see \cite{DJMP2,DJMP3,PV}.

Bipartite $(3,3,-4)$-graphs may be irregular. Figure \ref{fig:Bipartite(3,3,-4)} depicts all such graphs, which were obtained by using the program {\it geng} from the package {\it nauty} written by McKay~\cite{Mck06}. The unique bipartite $(3,4,-4)$-graph is shown in Fig.~\ref{fig:Bipartite(3,4,-4)}.

\begin{figure}[phtb]
\begin{center}
\includegraphics[scale=.75]{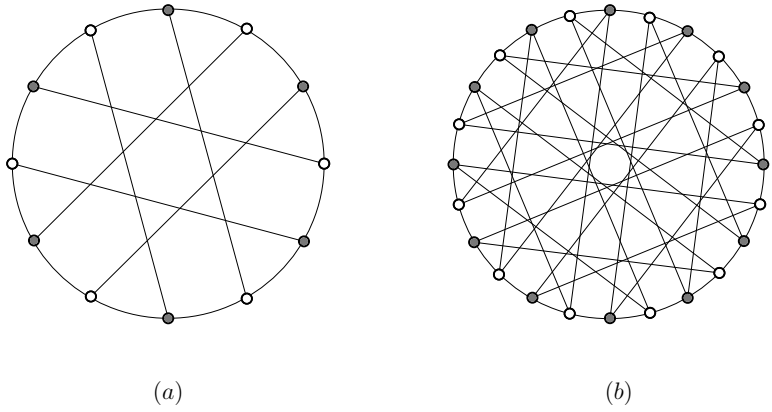}
\caption{($a$) the unique bipartite $(3,3,-2)$-graph and ($b$) the unique bipartite $(4,3,-2)$-graph.}
\label{fig:Bipartite(d,3,-2)}
\end{center}
\end{figure}

\begin{figure}[phtb]
\begin{center}
\makebox[\textwidth][c]{\includegraphics[scale=.75]{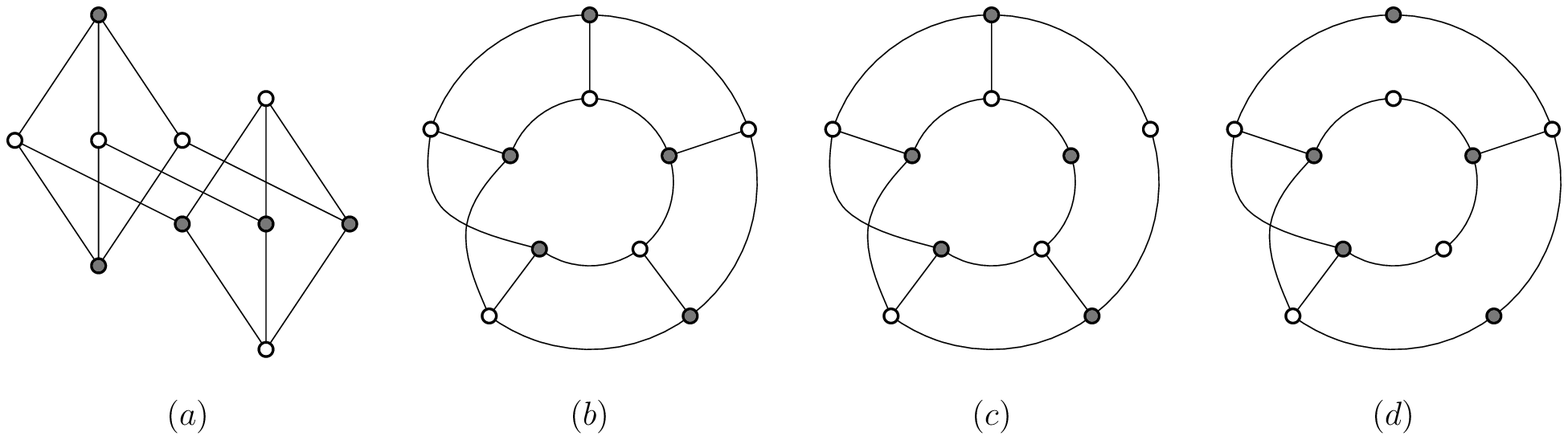}}
\caption{All the bipartite $(3,3,-4)$-graphs.}
\label{fig:Bipartite(3,3,-4)}
\end{center}
\end{figure}

\begin{figure}[phtb]
\begin{center}
\makebox[\textwidth][c]{\includegraphics[scale=.75]{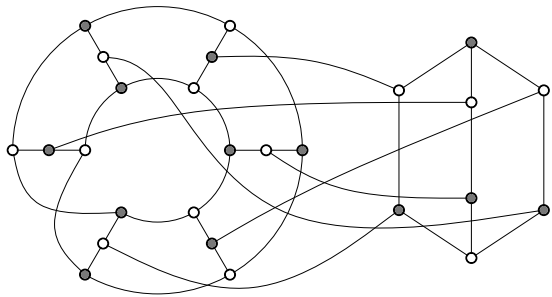}}
\caption{The unique bipartite $(3,4,-4)$-graph.}
\label{fig:Bipartite(3,4,-4)}
\end{center}
\end{figure}

\begin{figure}[phtb]
\begin{center}
\makebox[\textwidth][c]{\includegraphics[scale=.75]{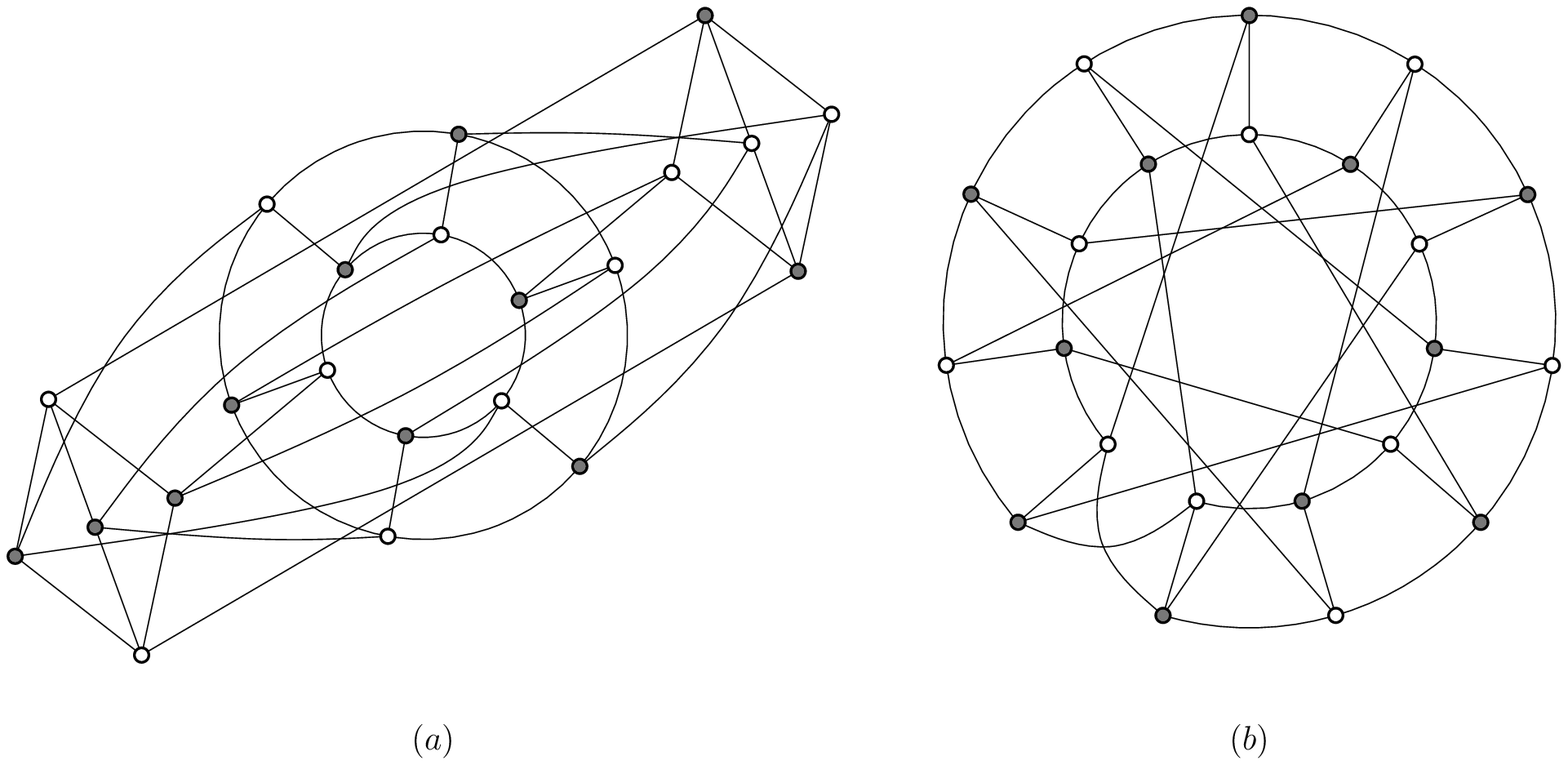}}
\caption{All the bipartite $(4,3,-4)$-graphs.}
\label{fig:Bipartite(4,3,-4)}
\end{center}
\end{figure}

\begin{figure}[phtb]
\begin{center}
\makebox[\textwidth][c]{\includegraphics[scale=.75]{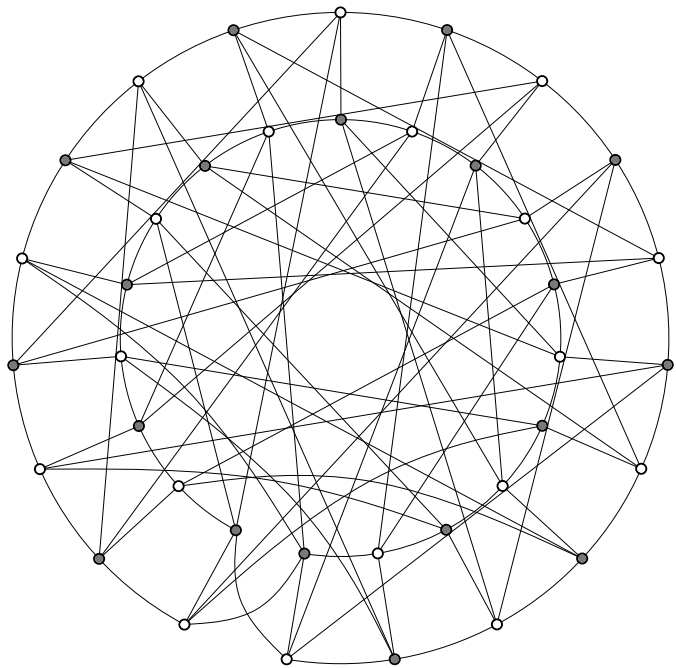}}
\caption{The only known bipartite $(5,3,-4)$-graph.}
\label{fig:Bipartite(5,3,-4)}
\end{center}
\end{figure}

All the bipartite $(4,3,-4)$-graphs are depicted in Fig.~\ref{fig:Bipartite(4,3,-4)}. These graphs were obtained computationally by Meringer \cite{Mer99} using the program {\it genreg}. An alternative description of the graph in Fig.~\ref{fig:Bipartite(4,3,-4)} ($b$) was communicated to the second author by Charles Delorme: take $\mathbb{Z}/22\mathbb{Z}$ as the vertex set of the graph, and for each even $x$, add the edges $\{x,x+1\}$, $\{x,x-1\}$, $\{x,x+7\}$ and $\{x,x+11\}$.

The only known bipartite $(5,3,-4)$-graph is depicted in Fig.~\ref{fig:Bipartite(5,3,-4)}; this graph was independently found by Charles Delorme and by the first author. Charles Delorme described this graph as follows: take $\mathbb{Z}/38\mathbb{Z}$ as its
vertex set, and for each even $x$, add the edges $\{x,x-1\}$, $\{x,x+1\}$, $\{x,x+5\}$, $\{x,x+13\}$ and $\{x,x+23\}$.


\section{Preliminary Results}
\label{sec:PreResults}

From now on we use the symbol $d$ rather than $\D$ to denote the degree of a regular graph, as it is customary. Recall that, unless $d=3$ and $D=3$, a bipartite $(d,D,-4)$-graph with $\D\ge3$ and $D\ge 3$ must be regular. Therefore, when referring to a regular bipartite $(d,D,-4)$-graph with $\D\ge3$ and $D\ge 3$, we are actually referring to any bipartite $(d,D,-4)$-graph with $\D\ge3$ and $D\ge 3$ other than the ones exhibited in Fig.~\ref{fig:Bipartite(3,3,-4)} $(c)$ and $(d)$.

In a bipartite $(d,D,-4)$-graph we call a cycle of length at most $2D-2$ a {\it short cycle}.

\begin{proposition}
\label{prop:GirthBipartite(d,D,-4)} The girth of a regular bipartite $(d,D,-4)$-graph $\G$ with $d\ge3$ and $D\ge 3$ is $2D-2$. Furthermore, any vertex $x$ of $\G$ lies on the short cycles specified below and no other short cycle, and we have the following cases:

{\bf $x$ is contained in {\bf exactly three $(2D-2)$-cycles}}. Then
\begin{itemize}
\item[{\rm ($i$)}] $x$ is a branch vertex of one $\Theta_{D-1}$, or
\end{itemize}
{\bf $x$ is contained in {\bf two $(2D-2)$-cycles}}. Then

\begin{itemize}
\item[{\rm ($ii$)}] $x$ lies on {\bf exactly two $(2D-2)$-cycles}, whose
intersection is a $\ell$-path with $\ell\in \{0,\ldots,D-1\}$.
\end{itemize}
\end{proposition}
Each case is considered as a type. For instance, a vertex satisfying
case $(i)$ is called a vertex of Type $(i)$. Note that, if $x$ is of Type $(ii)$ and $\ell=D-1$, the two short cycles containing $x$ constitute a $\Theta_{D-1}$.

\pr Let $xy$ be an edge of $\G$. Let us use the standard decomposition for a bipartite graph of even girth with respect to the edge $xy$ \cite{BI80}. For $0\le i\le D-1$, the sets $X_i$ and $Y_i$ are defined as follows:

\begin{center}
\begin{tabular}{lll}
$X_i$&=&$\{z\in V({\G})|d(x,z)=i, d(y,z)=i+1\}$\\
$Y_i$&=&$\{z\in V({\G})|d(y,z)=i, d(x,z)=i+1\}$.
\end{tabular}
\end{center}

The decomposition of $\G$ into the sets $X_i$ and $Y_i$ is called the {\it
standard decomposition for a graph of even girth with respect to the edge
$xy$}. Since $\G$ is bipartite, its girth is even and $X_i\cap
Y_j=\emptyset$ for $0\le i,j\le D-1$.

{\bf Claim 1} $\g(\G)= 2D-2$.

{\bf Proof of Claim 1.} Since the assertion is trivial for $D=3$, we suppose that $\g(\G)\le 2D-4$ for $D\ge4$. Assume that the edge $xy$ lies on a cycle of length $\g(\G)$. Then, $|X_i|=|Y_i|= (d-1)^i$ for $1 \leq i \le \frac{\g(\G)}{2}-1$, and
 \begin{eqnarray*}
 |X_{D-2}|&\leq& \left((d -1)^{D-3}-1\right)(d-1)+d-2=(d-1)^{D-2}-1\\
  |Y_{D-2}|&\leq& \left((d -1)^{D-3}-1\right)(d-1)+d-2=(d-1)^{D-2}-1\\
   |X_{D-1}|&\leq& \left((d-1)^{D-2}-1\right)(d-1)\\
  |Y_{D-1}|&\leq& \left((d-1)^{D-2}-1\right)(d-1).
 \end{eqnarray*}
Therefore,
\begin{eqnarray*}
|V(\G)|=\sum_{i=0}^{D-1} |X_i|+ \sum_{i=0}^{D-1} |Y_i| & \le & 2\left(1 +
(d-1)+(d -1)^2 + \dots +(d -1)^{D-3}\right) +{}\nonumber\\
&& +2(d-1)^{D-2}-2+2(d-1)^{D-1}-2(d-1)=\nonumber\\
& = & 2\left(1 +(d-1)+(d -1)^2 + \dots +(d -1)^{D-1}\right)-2(d-1)-2=\nonumber\\
&=& \M^b(d,D)-2d\nonumber,
\end{eqnarray*}

which is a contradiction. Hence, $\g(\G)\ge 2D-2$. If $\g(\G)= 2D$ then the order of $\G$ would be at least $\M^b(d,D)$ \cite{Big93}. Thus, $\g(\G)= 2D-2$ and the claim follows.\EndProof

We now proceed to prove the second part of the proposition.

For a given vertex $x$, we use again the standard decomposition for a bipartite graph with respect to an edge $xy$ in $\G$. Suppose that there are at least three edges joining vertices at $X_{D-2}$ to vertices at $Y_{D-2}$; that is, $|E(X_{D-2},Y_{D-2})|\ge3$ . In such case

\begin{center}
$|X_{D-1}|\leq (d -1)^{D-1}-3$,

$|Y_{D-1}|\leq (d -1)^{D-1}-3$,
\end{center}

 and therefore
\begin{align*}
|V(\G)|=\sum_{i=0}^{D-1} |X_i|+ \sum_{i=0}^{D-1} |Y_i|  \le {} & 2\left(1 +
(d-1)+(d -1)^2 + \dots +(d -1)^{D-2}\right) +{}\nonumber\\
& +2\left((d-1)^{D-1}-3\right)=\nonumber\\
= {} & 2\left(1 +(d-1)+(d -1)^2 + \dots +(d -1)^{D-1}\right)-6=\nonumber\\
= {} & M^b_{d,D}-6\nonumber,
\end{align*}

which is a contradiction. Consequently, $0\le|E(X_{D-2},Y_{D-2})|\le 2$.

Suppose that $|E(X_{D-2},Y_{D-2})|=2$. If the two edges are both incident to a common vertex of $Y_{D-2}$ then $x$ is of Type $(i)$, otherwise $x$ is of Type $(ii)$.

If instead $|E(X_{D-2},Y_{D-2})|=1$ then $|E(X_{D-2}, X_{D-1})|=|E(Y_{D-2}, Y_{D-1})|=(d-1)^{D-1}-1$. Since $|X_{D-1}|=|Y_{D-1}|=(d-1)^{D-1}-2$, there is a vertex $u\in X_{D-1}$ such that $|E(u,X_{D-2})|=2$. Therefore, it follows $(ii)$.

Finally, if $|E(X_{D-2},Y_{D-2})|=0$ then both types may occur. Indeed, if there is a vertex $u\in X_{D-1}$ such that $|E(u,X_{D-2})|=3$ then $x$ is of Type $(i)$ (this case can only occur if $d\ge4$), otherwise there must exist two vertices $u,v\in X_{D-1}$ such that $|E(u,X_{D-2})|=|E(v,X_{D-2})|=2$, in which case $x$ is of Type $(ii)$.  This completes the proof of the proposition.\EndProof

We continue with the following observation, which will be implicitly used throughout the paper:

\begin{observation}
\label{obs:Bipartite1} Let $\G =(V_1\cup V_2, E)$ (the sets $V_1$ and $V_2$ are called {\rm partite sets}) be any bipartite
graph of even $(odd)$ finite diameter $D$. The distance between a
vertex $u\in V_1$ and any vertex $v\in V_2$ $(w\in V_1)$ is at most
$D-1$.
\end{observation}

In virtue of Proposition \ref{prop:GirthBipartite(d,D,-4)}, we define the following concepts:

If two short cycles $C^1$ and $C^2$ are non-disjoint we say that $C^1$ and $C^2$ are \emph{neighbours}.

For a vertex $x$ lying on a short cycle $C$, we denote by $\rep^C(x)$ the vertex $x'$ in $C$ such that $d(x,x')=D-1$. We say $x'$ is the \emph{repeat} of $x$ in $C$ and vice versa, or simply that $x$ and $x'$ are \emph{repeats} in $C$. Alternatively, and more generally, we say that $x'$  is a repeat of $x$ with {\it multiplicity} $\m_x(x')$ $(1 \leq \m_x(x')\le2)$ if there are exactly $\m_x(x')+1$ different paths of length $D-1$ from $x$ to $x'$. Proposition \ref{prop:GirthBipartite(d,D,-4)} tells us that a vertex in $\G$ may have a repeat of multiplicity 2. Accordingly, we denote by $Rep(x)$ the multiset of the repeats  of a vertex $x$ in $\G$.

The concept of repeat can be easily extended to paths. For a path $P=x-y$ of length at most $D-2$ contained in a short cycle $C$, we denote by $\rep^C(P)$ the path $P'\subset C$ defined as $\rep^C(x)-\rep^C(y)$. We say that $P'$ is the \emph{repeat} of $P$ in $C$ and vice versa, or simply that $P$ and $P'$ are \emph{repeats} in $C$.

Often our arguments revolve around the identification of the elements in the set $S_x$ of short cycles containing a given vertex $x$; we call this process \emph{saturating} the vertex $x$. A vertex $x$ is called \emph{saturated} if the elements in $S_x$ have been completely identified. The following lemma will help us in this cycle identification process.


\begin{lemma}[Saturating Lemma]
\label{lemm:SaturatingLemma} Let $\C$ be a $(2D-2)$-cycle in a regular bipartite $(d,D,-4)$-graph $\G$ with $d \ge 3$ and $D \ge 3$, and $\alpha,\alpha'$ two vertices in $\C$ such that $\alpha'=\rep^{\C}(\alpha)$. Let $\gamma$ be a neighbour of $\alpha$ not contained in $\C$, and $\mu_1,\mu_2,\ldots,\mu_{d-2}$ the neighbours of $\alpha'$ not contained in $\C$. Suppose there is no short cycle in $\G$ containing the edge $\alpha\sim\gamma$ and intersecting $\C$ at a path of length greater than $D-3$.

Then, in $\G$ there exist  a vertex $\mu\in\{\mu_1,\mu_2,\ldots,\mu_{d-2}\}$ and a short cycle $\C^1$ such that $\gamma$ and $\mu$ are repeats in $\C^1$, and $\C\cap\C^1=\emptyset$.
\end{lemma}

\pr  Let $\alpha'_1,\alpha'_2$ be the neighbours of $\alpha'$ contained in $\C$.

\begin{figure}[!ht]
\begin{center}
\makebox[\textwidth][c]{\includegraphics[scale=1]{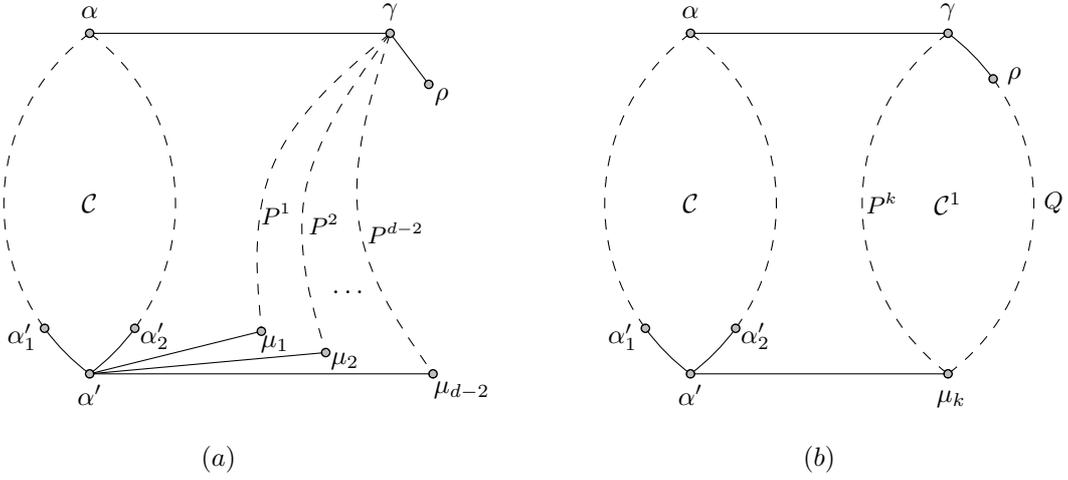}}
\caption{Auxiliary figure for Lemma \ref{lemm:SaturatingLemma}}
\label{Single Saturating Figure}
\end{center}
\end{figure}

For $1\le i\le d-2$, consider the path $P^i=\gamma-\mu_i$. As $\g(\G)=2D-2$, $P^i$ cannot go through $\alpha'_1$ or $\alpha'_2$. If $P^i$ went through certain $\mu_j$ $(j\ne i)$, then a cycle $\gamma P^i\mu_j\alpha'\C\alpha\gamma$ would either have length smaller than $2D-2$ or be a short cycle intersecting $\C$ at a $(D-1)$-path. Consequently, $P^i$ must go through one of the neighbours of $\mu_i$ other than $\alpha'$, and must be a $(D-1)$-path; see Fig.~\ref{Single Saturating Figure} $(a)$. In addition, $V(P^i\cap\C)=\emptyset$.

Let $\rho$ be one of the neighbours of $\gamma$ other than $\alpha$, not contained in any of the paths $P^i$ (there is at least one of such vertices). Consider a path $Q=\rho-\alpha'$. If $Q$ went through $\alpha'_1$, then the closed walk $\rho Q\alpha'_1\C\alpha\gamma\rho$ would either contain a cycle of length smaller than $2D-2$ or be a short cycle intersecting $\C$ at a $(D-2)$-path. Consequently, $Q$ must go through a certain $\mu_k$ $(1\le k\le d-2)$ and $V(Q\cap\C)=\{\alpha'\}$ (Fig.~\ref{Single Saturating Figure} $(b)$). Note that $Q$ must be a $(D-1)$-path, and that $V(Q\cap P^k)=\{\mu_k\}$; otherwise there would be a cycle in $\G$ of length smaller than $2D-2$.

Thus, we have obtained a short cycle $\C^1=\gamma\rho Q\mu_kP^k\gamma$ such that $\gamma$ and $\mu_k$ are repeats in $\C^1$, and $\C\cap\C^1=\emptyset$. By setting $\mu=\mu_k$ the lemma follows. \EndProof

\begin{corollary}
\label{cor:PreNeighborhoodLemma} Let $\alpha,\gamma$ be vertices in a regular bipartite $(d,D,-4)$-graph $\G$ with $d \ge 3$ and $D \ge 3$ such that $\gamma\in N(\alpha)$. Then, for every $\alpha'\in Rep(\alpha)$ it follows that $N(\alpha')$ contains a repeat of $\gamma$.
\end{corollary}

\pr Let $\C$ be a short cycle containing $\alpha$ and $\alpha'$. If the vertex $\gamma$ is contained in $\C$ or the edge $\alpha\gamma$ belongs to a short cycle in $\G$ intersecting $\C$ at a path of length $D-2$ or $D-1$, then the corollary trivially follows. If we instead assume that $\gamma\not\in\C$ and there is no short cycle in $\G$ containing the edge $\alpha\gamma$ and intersecting $\C$ at a path of length greater than $D-3$, then the corollary follows from the Saturating Lemma. \EndProof


\subsection{Repeats of Cycles}
In this section we extend the concept of repeat to short cycles; see the Repeat Cycle Lemma.

\begin{lemma}[Repeat Cycle Lemma]
\label{lemm:RepeatCycleLemma}
Let $C$ be a short cycle in a regular bipartite $(d,D,-4)$-graph $\G$ with $d \ge 3$ and $D \ge 3$, $\{C^1,C^2,\ldots ,C^k\}$ the set of neighbours of $C$, and  $I_i=C^i\cap C$ for $1\le i\le k$. Suppose at least one $I_j$, for $j\in\{1,\ldots, k\}$, is a path of length smaller than $D-2$.  Then there is an additional short cycle $C'$ in $\G$ intersecting $C^i$ at $I'_i=\rep^{C^i}(I_i)$, where $1\le i\le k$.
\end{lemma}

\begin{figure}[!ht]
\begin{center}
\makebox[\textwidth][c]{\includegraphics[scale=.9]{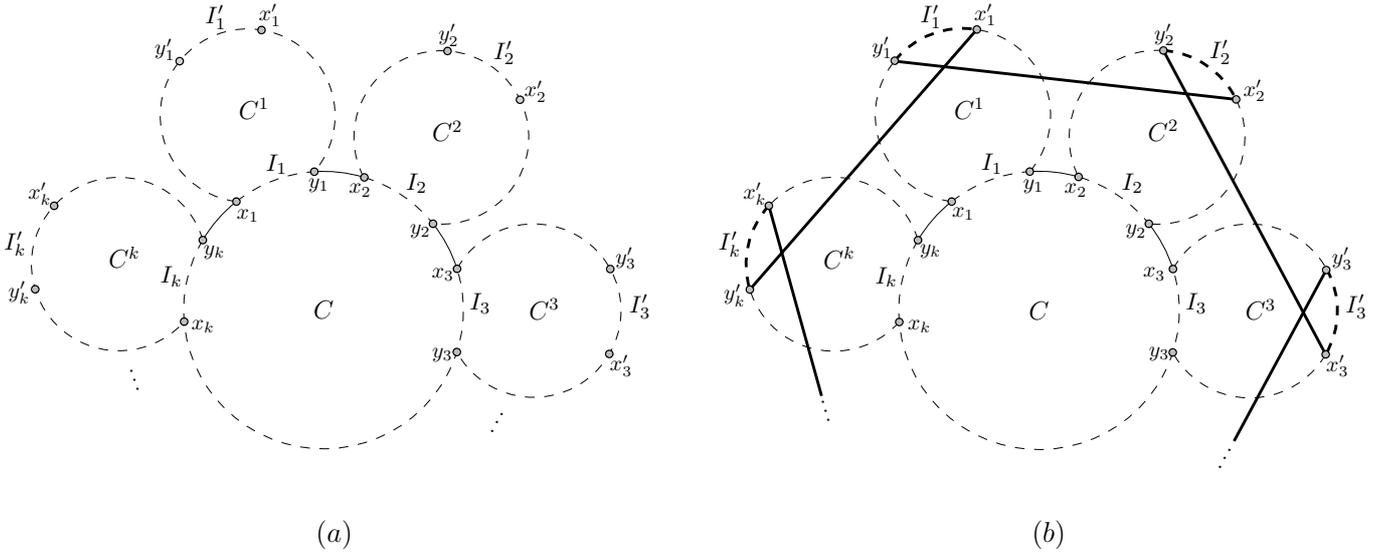}}
\caption{Auxiliary figure for Lemma \ref{lemm:RepeatCycleLemma}}
\label{Repeat Cycle Figure}
\end{center}
\end{figure}

\pr Observe that, according to our premises and Proposition \ref{prop:GirthBipartite(d,D,-4)}, $k\ge3$ and $I_i \cap I_j = \emptyset$ for $1\le i < j \le k$. We assume the denotation of the neighbours $C^1,C^2,\ldots,C^k$ of $C$ and the corresponding intersection paths $I_1=x_1-y_1,I_2=x_2-y_2,\ldots,I_k=x_k-y_k$ is such that $C=x_1I_1y_1x_2I_2y_2\ldots x_kI_ky_kx_1$. For $1\le i\le k$, we also denote the endvertices of $I'_i$ by $x'_i$ and $y'_i$, where $x'_i=\rep^{C^i}(x_i)$ and $y'_i=\rep^{C^i}(y_i)$ (see Fig.~\ref{Repeat Cycle Figure} ($a$)).

For $1\le i \le k$, consider the cycles $C^i$ and $C^{(i \bmod{k})+1}$.

First suppose that $I_i$ is a path of length smaller than $D-2$. Since $y_i$ is saturated, there cannot be a short cycle in $\G$, other than $C$, containing the edge $y_i\sim x_{(i \bmod{k})+1}$. Since $I_i$ is a path of length smaller than $D-2$, we apply the Saturating Lemma (mapping $C^i$ to $\C$, $y_i$ to $\alpha$, $y'_i$ to $\al'$ and $x_{(i \bmod{k})+1}$ to $\gamma$) and obtain an additional short cycle $\C^1$ in $\G$ such that $x_{(i \bmod{k})+1}$ is a repeat  in $\C^1$ of a neighbour $v\not\in C^i$ of $y'_i$, and $\C^1\cap C_i=\emptyset$. Since $x_{(i \bmod{k})+1}$ is saturated, we have that necessarily $\C^1=C^{(i \bmod{k})+1}$, which in turn implies $v=x'_{(i \bmod{k})+1}$. In other words, it follows that $y'_i\sim x'_{(i \bmod{k})+1}\in E(\G)$.

If instead $I_i$ is a $(D-2)$-path then $I_{(i \bmod{k})+1}$ must be a path of length smaller than $D-2$. Therefore, we can apply the above reasoning and deduce that $x'_{(i \bmod k)+1}\sim y'_i \in E(\G)$.

This way we obtain a subgraph $\Upsilon= \bigcup_{i=1}^{k} \big( I'_i \cup y'_i\sim x'_{(i \bmod k)+1} \big)=x'_1I'_1y'_1x'_2I'_2y'_2\ldots x'_kI'_ky'_kx'_1$
intersecting $C^i$ at $I'_i$ for $1\le i\le k$ (see Fig.~\ref{Repeat Cycle Figure} ($b$), where part of the subgraph $\Upsilon$ is highlighted in bold).

We next show that $\Upsilon$ must be indeed a cycle.

{\bf Claim 1.} $\Upsilon$ is a $(2D-2)$-cycle.

{\bf Proof of Claim 1.}

If the paths $I'_i$ are pairwise disjoint then $\Upsilon$ is obviously a $(2D-2)$-cycle.

Suppose the paths $I'_i$ are not pairwise disjoint; then $|V(\Upsilon)|<2D-2$ and, according to Proposition \ref{prop:GirthBipartite(d,D,-4)}, $\Upsilon$ contains no cycle. Since $\Upsilon$ is clearly connected it is therefore a tree.

Let $z \in C^\ell$ be an arbitrary leaf in $\Upsilon$. If the repeat path $I'_\ell=x'_\ell-y'_\ell$ had length greater than 0, then $z$ would have at least two neighbours in $\Upsilon$. Therefore, $I_\ell=C \cap C^\ell$ contains exactly one vertex, and thus, $x_\ell = y_\ell$ and $z = x'_\ell = y'_\ell$.

Recall we do addition modulo $k$ on the subscripts of the vertices and the superscripts of the cycles.

\begin{figure}[!ht]
\begin{center}
\makebox[\textwidth][c]{\includegraphics[scale=1]{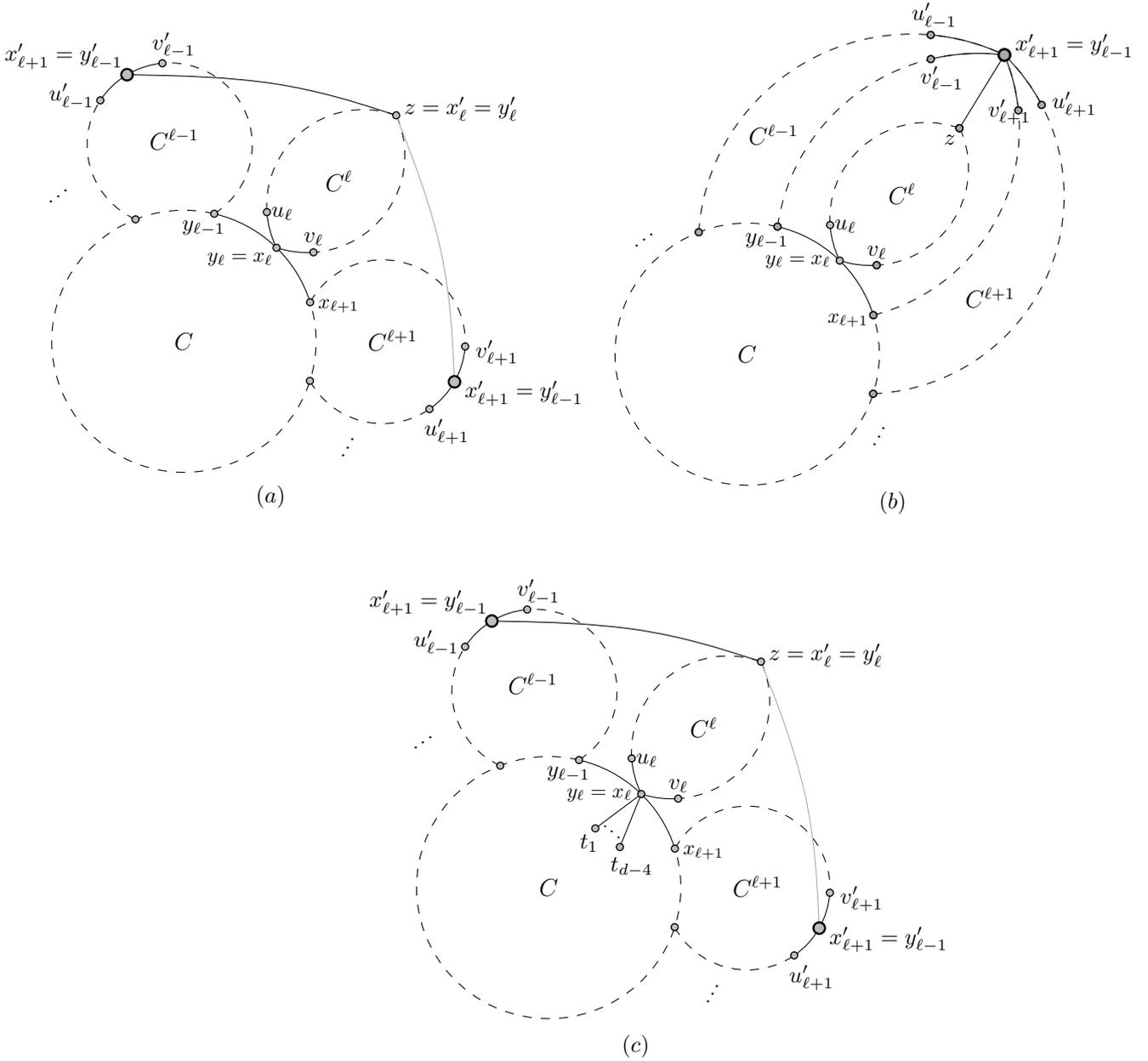}}
\caption{Auxiliary figure for Claim 1 of Lemma \ref{lemm:RepeatCycleLemma}.}
\label{fig:RepeatCyclePatch}
\end{center}
\end{figure}

Since $x'_\ell \sim y'_{\ell-1}$ and $x'_\ell \sim x'_{\ell+1}$ are edges in $\Upsilon$, it holds that $y'_{\ell-1}$ and $x'_{\ell+1}$ denote the same vertex. Let $u'_{\ell-1}, v'_{\ell-1}$ be the neighbours of $y'_{\ell-1}$ in $C^{\ell-1}$; $u'_{\ell+1}, v'_{\ell+1}$ the neighbours of $x'_{\ell+1}$ in $C^{\ell+1}$; and $u_\ell, v_\ell$ the neighbours of $x_\ell$ in $C^\ell$. We have that $V(C^{\ell-1} \cap C^{\ell+1}) = \{y'_{\ell-1}\}$, otherwise there would be a third short cycle in $\G$ containing $x_\ell$. In particular, the vertices in $\{u'_{\ell-1}, v'_{\ell-1}, u'_{\ell+1}, v'_{\ell+1}, x'_\ell \}$ are pairwise distinct and $d\ge5$. See Fig.~\ref{fig:RepeatCyclePatch} ($a$) and ($b$) for two drawings of this situation.

Let $t_1,t_2,\ldots,t_{d-4}$ denote the vertices in $N(x_\ell)-\{y_{\ell-1}, x_{\ell+1}, u_\ell, v_\ell \}$; see Fig.~\ref{fig:RepeatCyclePatch} ($c$). Consider a path $Q^i=t_i-y'_{\ell-1}$. Recall that $Q^i$ has length at most $D-1$. Since $x_\ell$ cannot be contained in a further short cycle, $Q^i$ must be a $(D-1)$-path and go through a neighbour of $y'_{\ell-1}$ not contained in $\{u'_{\ell-1}, v'_{\ell-1}, u'_{\ell+1}, v'_{\ell+1}, x'_\ell \}$. Therefore, we have that $d\ge6$ and, by the pigeonhole principle, that there are two paths $Q^r$ and $Q^s$ containing a common neighbour of $y'_{\ell-1}$. This way, $x_\ell$ would be contained in a third short cycle, a contradiction.

As a result, we conclude that the repeat graph $\Upsilon$ of $C$ is indeed a $(2D-2)$-cycle $C'$ as claimed. This completes the proof of Claim 1, and thus, of the lemma. \EndProof

While not of primary interest, it is not difficult to prove now that the cycles $C^1,C^2,\ldots ,C^k$ in the previous lemma are pairwise disjoint.

We call the aforementioned cycle $C'$ the {\it repeat} of the cycle $C$ in $\G$, and denote it by $\rep(C)$. Next some simple consequences of the Repeat Cycle Lemma follow.

\begin{corollary}[Repeat Cycle Uniqueness]
\label{cor:RepeatCycleUniqueness}
If a short cycle $C$ has a repeat cycle $C'$ then $C'$ is unique.
\end{corollary}

\begin{corollary}[Repeat Cycle Symmetry]
\label{cor:RepeatCycleSymmetry}
If $C'=\rep(C)$ then $C=\rep(C')$.
\end{corollary}


\begin{corollary}
\label{cor:RepeatPath}
Let $C$ and $C^1$ be two short cycles in a regular bipartite $(d,D,-4)$-graph $\G$ with $d \ge 3$ and $D \ge 3$ which intersect at a path $I$ of length smaller than $D-2$, and let $I'=\rep^{C^{1}}(I)$. Then, the repeat cycle of $C$ intersects $C^1$ at $I'$.
\end{corollary}

\begin{corollary}[Handy Corollary]
\label{cor:HandyCorollary}
Let $\C$ be a short cycle in a regular bipartite $(d,D,-4)$-graph $\G$ with $d \ge 3$ and $D \ge 3$, and $x,x'$ repeat vertices in $\C$. Let $\C^1$ and $\C^2$ be the other short cycles containing $x$ and $x'$, respectively. Suppose that $I=\C^1\cap \C$ is a path of length smaller than $D-2$. Then, setting $y=\rep^{\C^1}(x)$ and $y'=\rep^{\C^2}(x')$, we have that $y$ and $y'$ are repeat vertices in the repeat cycle of $\C$.
\end{corollary}

\pr
We denote the $k$ neighbour cycles of $\C$ as  $E^1,E^2,\ldots E^k$ and their respective intersection paths with $\C$ as $I_1=x_1-y_1,I_2=x_2-y_2,\ldots,I_k=x_k-y_k$ in such way that $\C=x_1I_1y_1x_2I_2y_2\ldots x_kI_ky_kx_1$. For $1\le j\le k$, we also denote $I'_j=x'_j-y'_j$, where $x'_j=\rep^{E^j}(x_j)$ and $y'_j=\rep^{E^j}(y_j)$.

Obviously, for some $r,s$ $(1\le r,s\le k)$ we have that $\C^1=E^r$, $\C^2=E^s$, $x\in I_r$, $x'\in I_s$, $y\in I'_r,$ and $y'\in I'_s$. We may assume $r<s$. By the Repeat Cycle Lemma, the vertices $y$ and $y'$ belong to the repeat cycle $\C'$ of $\C$. Then the paths $xI_ry_rx_{r+1}I_{r+1}y_{r+1}\ldots x_{s-1}I_{s-1}y_{s-1}x_sI_sx'\subset \C$ and $yI'_ry'_rx'_{r+1}I'_{r+1}y'_{r+1}\ldots x'_{s-1}I'_{s-1}y'_{s-1}x'_sI'_sy'\subset \C'$ are both $(D-1)$-paths in $\G$, and the corollary follows. \EndProof

\begin{proposition}
\label{prop:CyclePartition}
The set $S(\G)$ of short cycles in a bipartite $(d,D,-4)$-graph $\G$ with $d \ge 3$ and $D \ge 3$ can be partitioned into sets $S_{D-1}(\G)$, $S_{D-2}(\G)$ and $S_{D-3}(\G)$, where

\begin{itemize}
\item[] $S_{D-1}(\G)$ is the set of short cycles whose intersections with neighbour cycles are $(D-1)$-paths;
\item[] $S_{D-2}(\G)$ is the set of short cycles whose intersections with neighbour cycles are $(D-2)$-paths; and
\item[] $S_{D-3}(\G)$ is the set of short cycles whose intersections with neighbour cycles are paths of length at most $D-3$.
\end{itemize}
\end{proposition}

\pr If $\G$ is one of the non-regular graphs in Fig.~\ref{fig:Bipartite(3,3,-4)} the result trivially follows. We then assume that $\G$ is regular.

Let $C$ be a short cycle in $\G$. If $C$ is contained in a $\Theta_{D-1}$ then, according to Proposition \ref{prop:GirthBipartite(d,D,-4)}, all the intersections of $C$ with its neighbour cycles are $(D-1)$-paths, in which case $C \in S_{D-1}(\G)$.

Now suppose that, for some short cycle $C^1$, $P_1=C\cap C^1$ is a path of length $D-2$. Note that all vertices in $P_1$ are saturated. Let $v$ be an arbitrary vertex in $P_1$, $v'=rep^C(v)$, and $C^2$ the short cycle other than $C$ containing $v'$. Suppose that $P_2=C\cap C^2$ is not a $(D-2)$-path. Then clearly $P_2$ cannot be a $(D-1)$-path, so it has length at most $D-3$. But according to Corollary \ref{cor:RepeatPath}, the cycle $rep(C^2)$ intersects $C$ at exactly $rep^C(P_2)$, a proper subpath of $P_1$. This implies that $rep(C^2)$ is a third short cycle containing the vertex $v$, a contradiction. Consequently, the intersections of $C$ with its (exactly two) neighbour cycles are $(D-2)$-paths, and $C \in S_{D-2}(\G)$.

Finally, if there is a short cycle intersecting $C$ at a path of length at most $D-3$ then, by the above reasoning, the intersections of $C$ with all of its neighbour cycles are paths of length at most $D-3$, and $C \in S_{D-3}(\G)$. \EndProof

The preceding result could be stated alternatively in term of vertices as follows:

\begin{proposition}
\label{prop:VertexPartition}
The set $V(\G)$ of vertices in a regular bipartite $(d,D,-4)$-graph $\G$ with $d \ge 3$ and $D \ge 3$ can be partitioned into sets $V_{D-1}(\G)$, $V_{D-2}(\G)$ and $V_{D-3}(\G)$, where

\begin{itemize}
\item[] $V_{D-1}(\G)$ is the set of vertices contained in cycles of  $S_{D-1}(\G)$;
\item[] $V_{D-2}(\G)$ is the set of vertices contained in cycles of  $S_{D-2}(\G)$;
\item[] $V_{D-3}(\G)$ is the set of vertices contained in cycles of  $S_{D-3}(\G)$;
\end{itemize}

and $S_{D-1}(\G)$, $S_{D-2}(\G)$, $S_{D-3}(\G)$ are defined as in Proposition \ref{prop:CyclePartition}. \EndProof
\end{proposition}

\section{Main results}
\label{sec:MainResults}
\subsection{Non-existence of subgraphs isomorphic to $\Theta_{D-1}$}


\begin{theorem}
\label{theo:Theta}
A bipartite $(d,D,-4)$-graph $\G$ with $d\ge 3$ and $D\ge 5$ does not contain a subgraph isomorphic to $\Theta_{D-1}$.
\end{theorem}

\begin{figure}[!ht]
\begin{center}
\makebox[\textwidth][c]{\includegraphics[scale=.9]{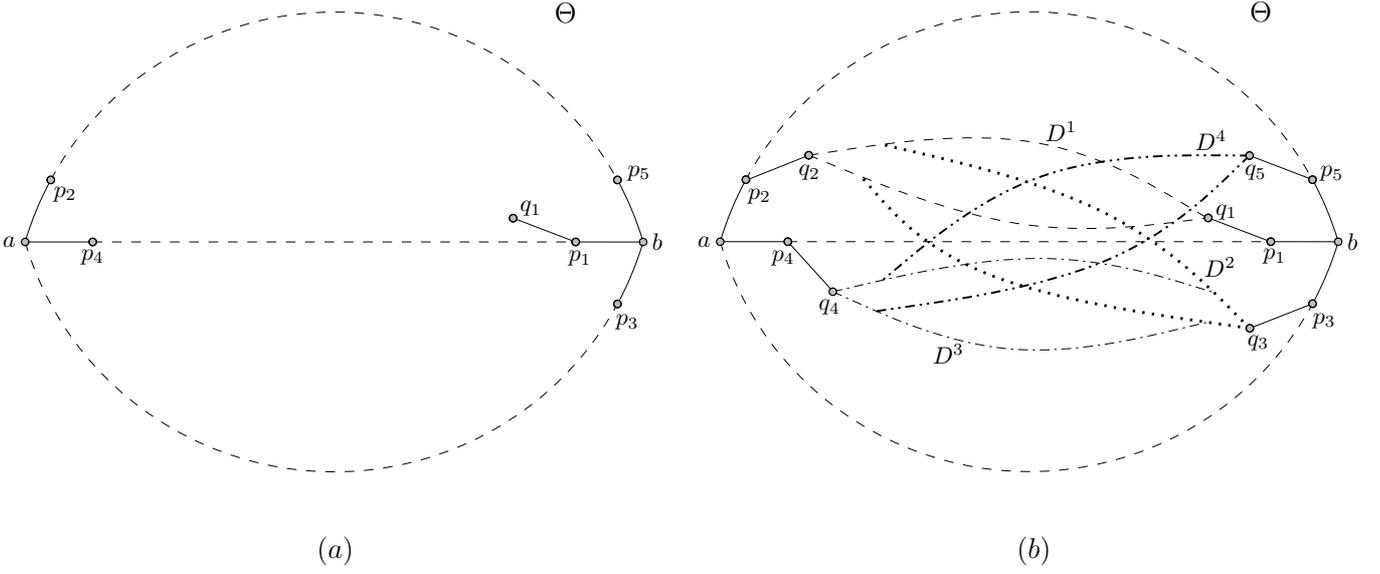}}
\caption{Auxiliary figure for Theorem \ref{theo:Theta}}
\label{Theta Figure}
\end{center}
\end{figure}

\pr
Suppose that $\G$ has a subgraph $\Theta$ isomorphic to $\Theta_{D-1}$, with branch vertices $a$ and $b$. Let $p_1,p_2,p_3,p_4$ and $p_5$ be as in Fig.~\ref{Theta Figure} ($a$), and let $q_1$ be one of the neighbours of $p_1$ not contained in $\Theta$.

Since all vertices of $\Theta$ are saturated, there cannot be a short cycle in $\G$ containing any of the incident edges of $p_1,p_2,p_3,p_4$ or $p_5$  which are not contained in $\Theta$. According to this and by applying the Saturating Lemma, there is an additional short cycle $D^1$ in $\G$ such that $q_1$ and one of the neighbours of $p_2$ not contained in $\Theta$ (say $q_2$) are repeats in $D^1$, and $D^1\cap \Theta=\emptyset$. Analogously,  in $\G$ there is an additional short cycle $D^2$ such that $q_2$ and one of the neighbours of $p_3$ not contained in $\Theta$ (say $q_3$) are repeats in $D^2$, and $D^2\cap \Theta=\emptyset$; an additional short cycle $D^3$ such that $q_3$ and one of the neighbours of $p_4$ not contained in $\Theta$ (say $q_4$) are repeats in $D^3$, and $D^3\cap \Theta=\emptyset$; and an additional short cycle $D^4$ such that $q_4$ and one of the neighbours of $p_5$ not contained in $\Theta$ (say $q_5$) are repeats in $D^4$, and $D^4\cap \Theta=\emptyset$. See  Fig.~\ref{Theta Figure} ($b$).

Note that $D^1\cap D^2$ is a path of length at most $2<D-2$; otherwise for some vertex $t\in D^1\cap D^2$ the closed walk $tD^1q_1p_1bp_3q_3D^2t$ would contain a cycle of length at most $2D-2$ to which the vertex $b$ would belong, a contradiction. For similar reasons, the intersection paths $D^2\cap D^3$ and $D^3\cap D^4$ all have length at most 2, with 2<$D-2$. We now apply the Handy Corollary. By mapping the cycle $D^2$ to $\C$, the vertex $q_2$ to $x$, the vertex $q_3$ to $x'$, the cycle $D^1$ to $\C^1$, the cycle $D^3$ to $\C^2$, the vertex $q_1$ to $y$ and the vertex $q_4$ to $y'$, we obtain that $q_1$ and $q_4$ are repeat vertices in the repeat cycle of $D^2$. Therefore, since $q_4\in D^4$,  it follows that $D^2$ and $D^4$ are repeat cycles and $q_1=q_5$. This way, there would be a cycle $q_1p_1bp_5q_5$ in $\G$ of length $4<2D-2$ (since $D\ge 5$), a contradiction to the fact that $\g(\G)=2D-2$.\EndProof

\begin{proposition}
\label{prop:NumberOfCycles}
The number $N_{2D-2}$ of short cycles in a bipartite $(d,D,-4)$-graph $\G$ with $d\ge 3$ and $D\ge 5$ is given by the expression $\frac{2\times\big(1+(d-1)+\ldots+{(d-1)}^{D-1}\big)-4}{D-1}$.
\end{proposition}

\pr
By Theorem \ref{theo:Theta}, $\G$ does not contain a subgraph isomorphic to $\Theta_{D-1}$. Then, according to Proposition \ref{prop:GirthBipartite(d,D,-4)}, every vertex of $\G$ is contained in exactly two short cycles. We then count the number $N_{2D-2}$ of short cycles of $\G$ . Since the order of $\G$ is $2\times\big(1+(d-1)+\ldots+{(d-1)}^{D-1}\big)-4$, we have that

\begin{center}
$N_{2D-2}=\frac{2\times\Big(2\times\big(1+(d-1)+\ldots+{(d-1)}^{D-1}\big)-4\Big)}{2D-2}=\frac{2\times\big(1+(d-1)+\ldots+{(d-1)}^{D-1}\big)-4}{D-1}$,
\end{center}

and the proposition follows.  \EndProof


\subsection{Non-existence results on bipartite $(d,D,-4)$-graphs}

Since the number of short cycles in a graph $\G$ must be an integer, the expression obtained for $N_{2D-2}$ in Proposition \ref{prop:NumberOfCycles} already suffices to prove the non-existence of bipartite $(d,D,-4)$-graphs for infinitely many pairs $(d,D)$.

Consider first the case in which $D-1=p^q$ is an odd prime power. Let $G=\{1,2,\ldots,p-1\}$ be the multiplicative group of the field $\mathbb{Z}/p\mathbb{Z}$, let $d-1 \not \equiv 0,1\pmod p$, and let $H$ be the cyclic subgroup of $G$ generated by $d-1$. We observe that the sum of the elements of $H$ is null $\pmod{p}$. Furthermore, since the order of $H$ divides the order of $G$, it must also divide $p^q-1=D-2$. Thus, we have
\begin{displaymath}
2\times\big(1+(d-1)+\ldots+{(d-1)}^{D-1}\big)-4\equiv
\left\{
\begin{array}{ll}
-2 \pmod p & \textrm{if $d-1\equiv 0,1\pmod p$,}\\
2(d-1)-2 \pmod p & \textrm{if $d-1\not\equiv 0,1\pmod p$.}
\end{array}
\right.
\end{displaymath}

Therefore, it immediately follows
\begin{corollary}
\label{cor:PrimePower}
There is no bipartite $(d,D,-4)$-graph with $d\ge3$ and $D\ge5$ such that $D-1$ is an odd prime power.\EndProof
\end{corollary}
More generally, if $p$ is an odd prime factor of $D-1$ and $D-1\equiv r\pmod{p-1}$, then
\begin{displaymath}
2\times\big(1+(d-1)+\ldots+{(d-1)}^{D-1}\big)-4\equiv
\left\{
\begin{array}{ll}
-2 \pmod p & \textrm{if $d-1\equiv 0,1\pmod p$,}\\
2\frac{(d-1)^{r+1}-1}{d-2}-4 \pmod p & \textrm{if $d-1\not\equiv 0,1\pmod p$;}
\end{array}
\right.
\end{displaymath}
\begin{corollary}
\label{cor:Residues01}
There is no bipartite $(d,D,-4)$-graph with $d\ge3$ and $D\ge6$ such that $d-1\equiv 0,1 \pmod p$, where $p$ is an odd prime factor of $D-1$.\EndProof
\end{corollary}

It is also possible to examine completely the case of some small odd prime factors of $D-1$. For example, it is not difficult to verify that, if $D-1=3k$ then $3$ does not divide $2\times\big(1+(d-1)+\ldots+{(d-1)}^{D-1}\big)-4$; thus,

\begin{corollary}
\label{cor:Dmod3}
There is no bipartite $(d,D,-4)$-graph with $d\ge3$ and $D\ge5$ such that $D-1\equiv 0 \pmod 3$.
\end{corollary}

Now we turn to structural arguments to obtain other non-existence results.
\begin{lemma}
\label{lemm:(D-2)-Intersection}
Any two non-disjoint short cycles in a bipartite $(d,D,-4)$-graph $\G$ with $d\ge 3$ and $D\ge 7$ intersect at a path of length smaller than $D-2$.
\end{lemma}

\pr
Since $\G$ does not contain a graph isomorphic to $\Theta_{D-1}$, it is only necessary to prove here that any two non-disjoint short cycles in $\G$ cannot intersect at a path of length $D-2$.

Suppose, by way of contradiction, that there are two short cycles $C^1$ and $C^2$ in $\G$ intersecting at a path $I_1$ of length $D-2$. According to Proposition \ref{prop:CyclePartition}, $C^2$ is intersected by exactly two short cycles, namely $C^1$ and $C^3$, at two independent $(D-2)$-paths. By repeatedly applying this reasoning and considering $\G$ is finite, we obtain a maximal length sequence $C^1,C^2,C^3,\ldots,C^{m}$ of pairwise distinct short cycles in $\G$ such that $C^i$ intersects $C^{i+1}$ at a path $I_i$ of length $D-2$ ($1\le i\le m-1$), and $C^i\cap C^j=\emptyset$ for any $i,j\in\{1,\ldots,m\}$ such that $2\le|i-j|\le m-2$.

\begin{figure}[!ht]
\begin{center}
\makebox[\textwidth][c]{\includegraphics[scale=.9]{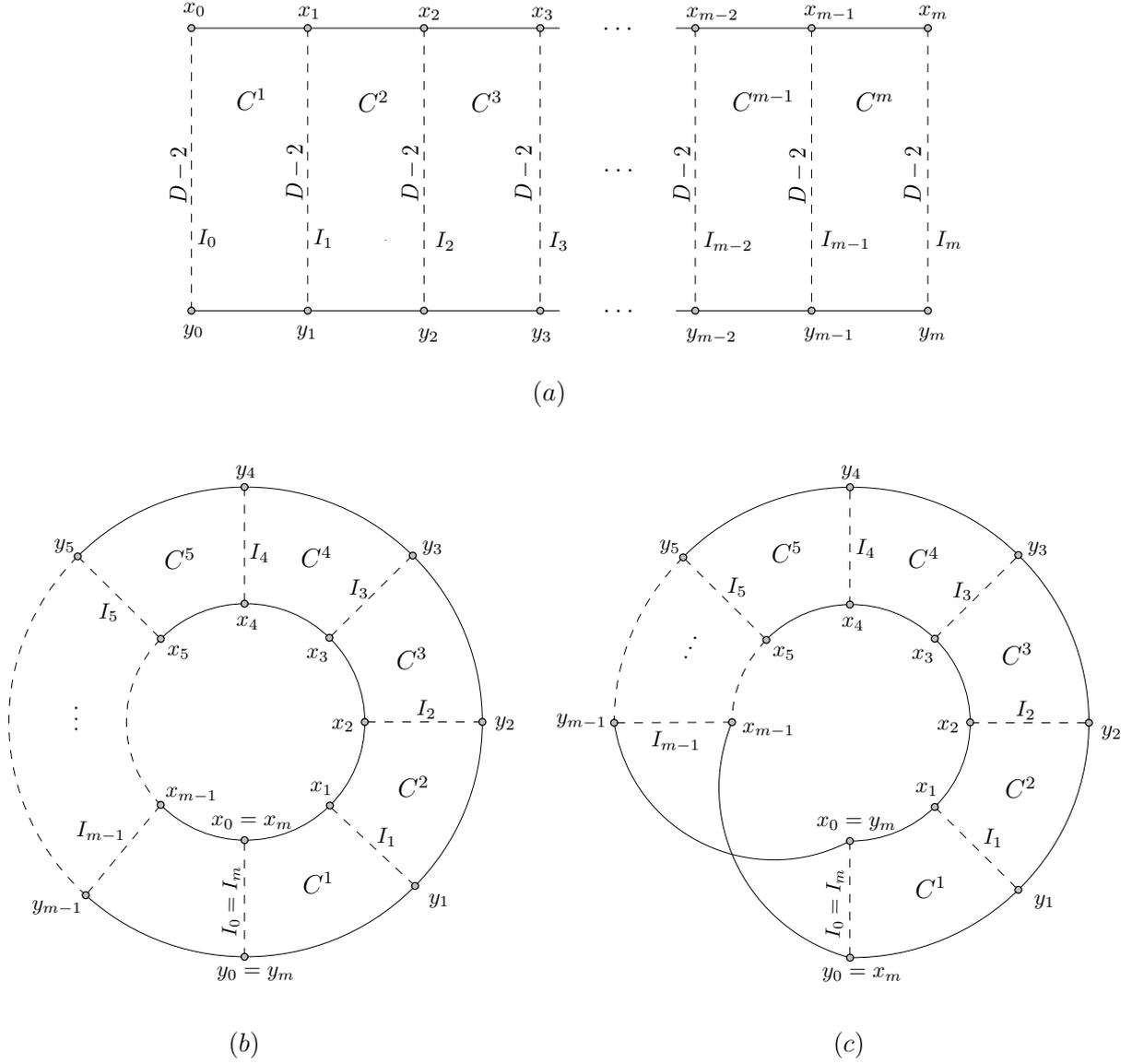}}
\caption{Auxiliary figure for Lemma \ref{lemm:(D-2)-Intersection}}
\label{Intersection Length1 Figure}
\end{center}
\end{figure}

Let us denote the paths $I_1=x_1-y_1,\ldots,$ $I_{m-1}=x_{m-1}-y_{m-1}$ in such way that, for $1\le i\le m-2$, $x_i\sim x_{i+1}$ and $y_{i}\sim y_{i+1}$ are edges in $\G$. Also, let $x_0\in N(x_{1})\cap (C^1-I_1)$, $y_0\in N(y_{1})\cap (C^1-I_1)$, $x_m\in N(x_{m-1})\cap (C^m-I_{m-1})$, and $y_m\in N(y_{m-1})\cap (C^m-I_{m-1})$; see Fig.~\ref{Intersection Length1 Figure} $(a)$. Set $I_0=x_0-y_0$ and $I_m=x_m-y_m$. Since the sequence $C^1,C^2,C^3,\ldots,C^{m}$ is maximal and all the vertices in $I_1,\ldots,I_{m-1}$ are saturated, it follows that $I_0=I_m$, and we have either $x_0=x_m$ and $y_0=y_m$ (Fig.~\ref{Intersection Length1 Figure} $(b)$), or $x_0=y_m$ and $y_0=x_m$ (Fig.~\ref{Intersection Length1 Figure} $(c)$).

If $x_0=x_m$ and $y_0=y_m$, then $m\ge 2D$; otherwise the cycle $x_1x_2\ldots x_mx_1$ would have length at most $2D-2$, contradicting the saturation of $x_1$. If conversely $x_0=y_m$ and $y_0=x_m$ then $m\ge D$; otherwise the cycle $x_1x_2\ldots x_my_1y_2\ldots y_mx_1$ containing $x_1$ would have length at most $2D-2$, a contradiction as well. For our purposes, it is enough to state $m\ge D\ge 7$ in any case.

\begin{figure}[!ht]
\begin{center}
\makebox[\textwidth][c]{\includegraphics[scale=.9]{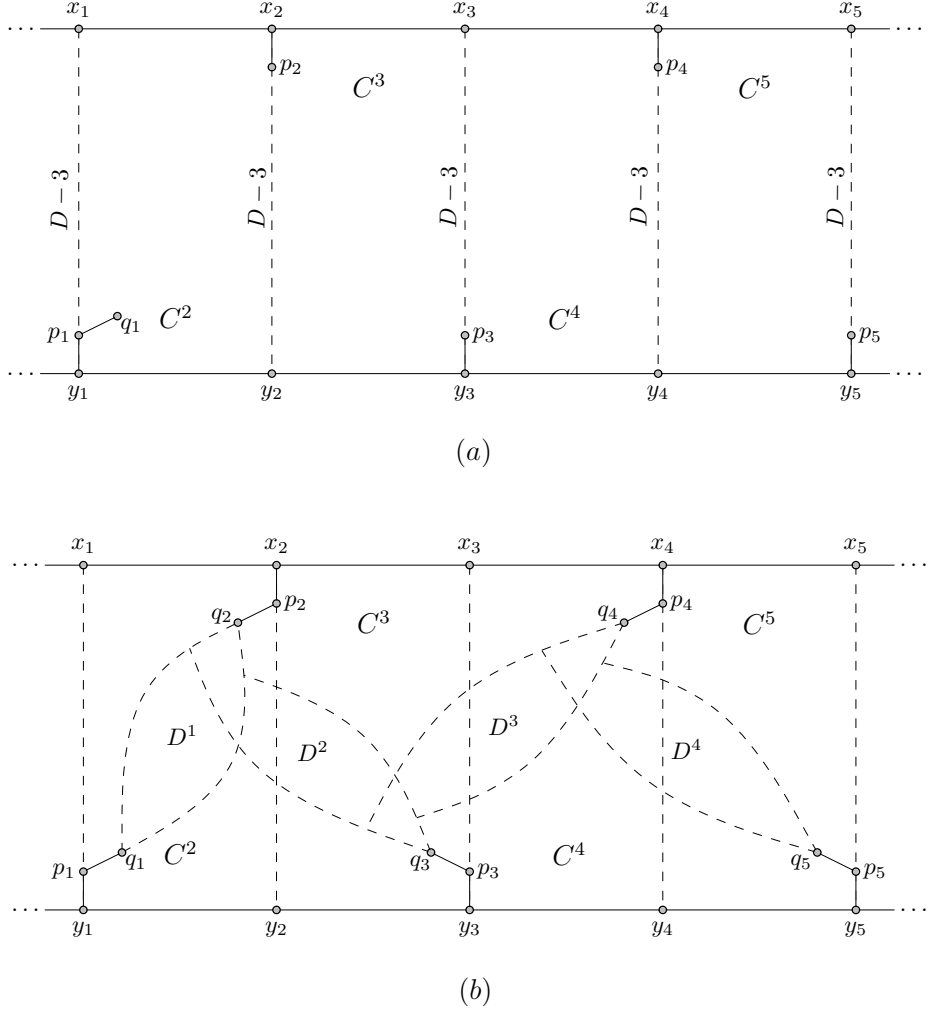}}
\caption{Auxiliary figure for Lemma \ref{lemm:(D-2)-Intersection}}
\label{Intersection Length2 Figure}
\end{center}
\end{figure}

Let $p_1$ be the neighbour of $y_1$ on $I_1$, and $p_{i+1}=\rep^{C^{i+1}}(p^i)$ for $1\le i\le 4$. Also, let $q_1$ be a neighbour of $p_1$ not contained in $I_1$; see Fig.~\ref{Intersection Length2 Figure} ($a$).

Since all vertices on $I_1$ are saturated, the edge $q_1\sim p_1$ cannot be contained in a further short cycle. We apply the Saturating Lemma (by mapping $C^{2}$ to $\C$, $p_1$ to $\alpha$, $p_{2}$ to $\al'$, and $q_1$ to $\gamma$), and obtain in $\G$ an additional short cycle $D^1$ such that $q_1$ and one of the neighbours of $p_{2}$ not contained in $I_{2}$ (say $q_{2}$) are repeats in $D^1$, and $D^1\cap C^2=\emptyset$. Analogously, for $2\le i\le 4$ we obtain an additional short cycle $D^i$ in $\G$ such that $q_i$ and a neighbour of $p_{i+1}$ not contained in $I_{i+1}$ (say $q_{i+1}$) are repeats in $D^i$, and $D^i\cap C^{i+1}=\emptyset$; see Fig.~\ref{Intersection Length2 Figure} ($b$).

For $i=1$ or $3$, $D^i\cap D^{i+1}$ cannot be a $(D-2)$-path; otherwise for some vertex $t_i\in D^i\cap D^{i+1}$, there would be a cycle $q_ip_iy_iy_{i+1}y_{i+2}p_{i+2}q_{i+2}D_{i+1}t_iD_iq_i$ of length at most $6+D-4+D-4$ (since $D-2\ge 5$), a contradiction to the fact that $p_i$ is saturated and $\g(\G)=2D-2$. Analogously, $D^2\cap D^3$ cannot be a $(D-2)$-path.

We now apply the Handy Corollary. By mapping  the cycles $D^2$ to $\C$, $D^1$ to $\C^1$ and $D^3$ to $\C^2$, and the vertices $q_2$ to $x$, $q_3$ to $x'$, $q_1$ to $y$, and $q_4$ to $y'$, it follows that the vertices $q_1$ and $q_4$ are repeat vertices in the repeat cycle of $D^2$. Since $q_4\in D^4$, $D^2$ and $D^4$ are repeat cycles and $q_5=q_1$. This way, we obtain a cycle $q_1p_1y_1y_2y_3y_4y_5p_5q_5$ in $\G$ of length $8<2D-2$, a contradiction.

This completes the proof of the lemma.\EndProof


\begin{theorem}
\label{theo:OddDiameter}
There are no bipartite $(d,D,-4)$-graphs for $d\ge 3$ and odd $D\ge 5$.
\end{theorem}

\pr The case $D=5$ can be easily discarded by using Proposition \ref{prop:NumberOfCycles}, so we assume $D\ge 7$.

Suppose there is a bipartite $(d,D,-4)$-graph $\G$ with $d\ge 3$ and odd $D\ge 7$. According to Lemma \ref{lemm:(D-2)-Intersection}, any two non-disjoint short cycles in $\G$ intersect at a path of length smaller than $D-2$, which means that every short cycle $C$ in $\G$ has a repeat cycle $C'$ (by the Repeat Cycle Lemma). Because of the uniqueness and symmetry of repeat cycles, the number $N_{2D-2}$ of short cycles in $\G$ must be even.

However, since $D$ is odd, the number $N_{2D-2}=\frac{2\times\big(1+(d-1)+\ldots+{(d-1)}^{D-1}\big)-4}{D-1}$ of short cycles in $\G$ is odd, a contradiction.\EndProof

Furthermore, using Theorem \ref{theo:OddDiameter} and Proposition \ref{prop:NumberOfCycles} we complete the catalogue of bipartite $(d,D,-4)$-graphs with $5\le D\le 187$. From Section \ref{sec:KnownGraphs} recall that, for $D\ge5$, the only bipartite $(2,D,-4)$-graph is the path of length $5$.

\begin{theorem}
\label{theo:5-187}
The path of length 5 is the only bipartite $(\D,D,-4)$-graph with $\D\ge2$ and $5\le D\le 187$.
\end{theorem}
\subsection{Non-existence of bipartite $(3,D,-4)$-graphs with $D\ge5$}

In this section we complete the catalogue of bipartite $(3,D,-4)$-graphs. Specifically, we prove the non-existence of bipartite $(3,D,-4)$-graphs with even $D\ge6$.

\begin{figure}[!ht]
\begin{center}
\makebox[\textwidth][c]{\includegraphics[scale=.9]{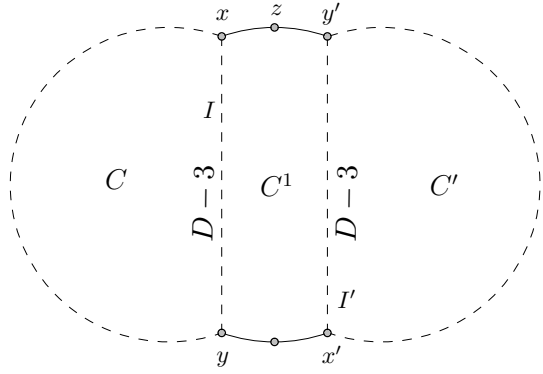}}
\caption{Auxiliary figure for Lemma \ref{lemm:D-3Int}}
\label{fig:D-3Int}
\end{center}
\end{figure}

\begin{lemma}\label{lemm:D-3Int} Any two non-disjoint short cycles in a bipartite $(3,D,-4)$-graph $\G$ with $d\ge 3$ and $D\ge 7$  intersect at a path of length smaller than $D-3$.
\end{lemma}
\pr By Lemma \ref{lemm:(D-2)-Intersection}, it is only necessary to prove here that any two short cycles $C$ and $C^1$ in $\G$ cannot intersect at a path $I=x-y$ of length $D-3$. We proceed by contradiction. Let $x'$ and $y'$ be the
repeat vertices of $x$ and $y$ in $C^1$, respectively. By Corollary \ref{cor:RepeatPath}, the repeat cycle $C'$ of $C$ intersects $C^1$ at $I'=x'-y'$ (the repeat path of $I$ in $C^1$); see Fig.~\ref{fig:D-3Int}. If we denote by $z$ the neighbour of $x$ on $C^1-C$, then we have that the other short cycle containing $z$ would also contain at least one of the vertices in $\{x,y'\}$, which contradicts the fact that $x$ and $y'$ are both saturated. \EndProof

\begin{theorem}
\label{theo:(3,D,-4)}
There are no bipartite $(3,D,-4)$-graphs with even $D\ge6$.
\end{theorem}

\pr Recall that Theorem \ref{theo:5-187} covers the case $D=6$.

Let $\G$ be a bipartite $(3,D,-4)$-graph with even $D\ge 8$, $C^0$ a short cycle in $\G$, and $x_0$, $x'_0$ two repeat vertices in $C^0$. Let $x_1$ and $x'_1$ be the neighbours of $x_0$ and $x'_0$, respectively, not contained in $C^0$. According to the Saturating Lemma, there is an additional short cycle $C^1$ containing $x_1$ and $x'_1$ such that $C^0\cap C^1=\emptyset$. Let $y_1$ be one of the neighbours of $x_1$ contained in $C^1$, and $y'_1=\rep^{C^1}(y_1)$. Denote by $x_2$ and $x'_2$ the neighbours of $y_1$ and $y'_1$, respectively, not contained in $C^1$. Again by  the Saturating Lemma, there is an additional short cycle $C^2$ such that $x'_2=\rep^{C^2}(x_2)$ and $C^1\cap C^2=\emptyset$. Since $d=3$, we may assume that the other short cycle $C$ containing $x_0$ also contains $x_1$ and a neighbour of $x_0$ in $C_0$ (say $y_0$). We first prove that $C^0\cap C=y_0x_0$.

\begin{figure}[!ht]
\begin{center}
\makebox[\textwidth][c]{\includegraphics[scale=.95]{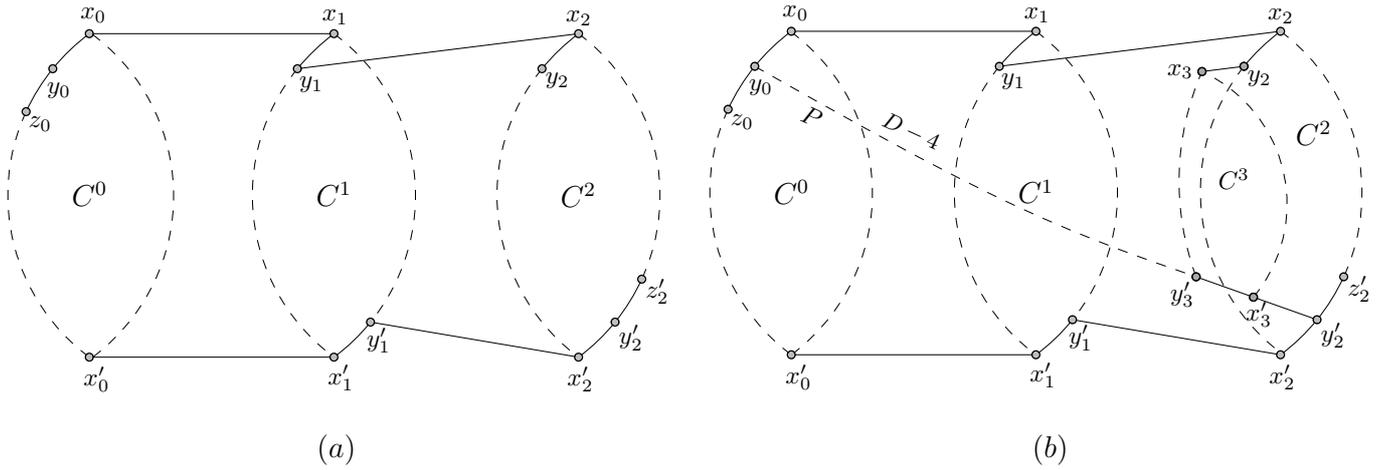}}
\caption{Auxiliary figure for Theorem \ref{theo:(3,D,-4)}.}
\label{fig:CubicEvenDiameter1}
\end{center}
\end{figure}

\begin{figure}[!ht]
\begin{center}
\makebox[\textwidth][c]{\includegraphics[scale=.9]{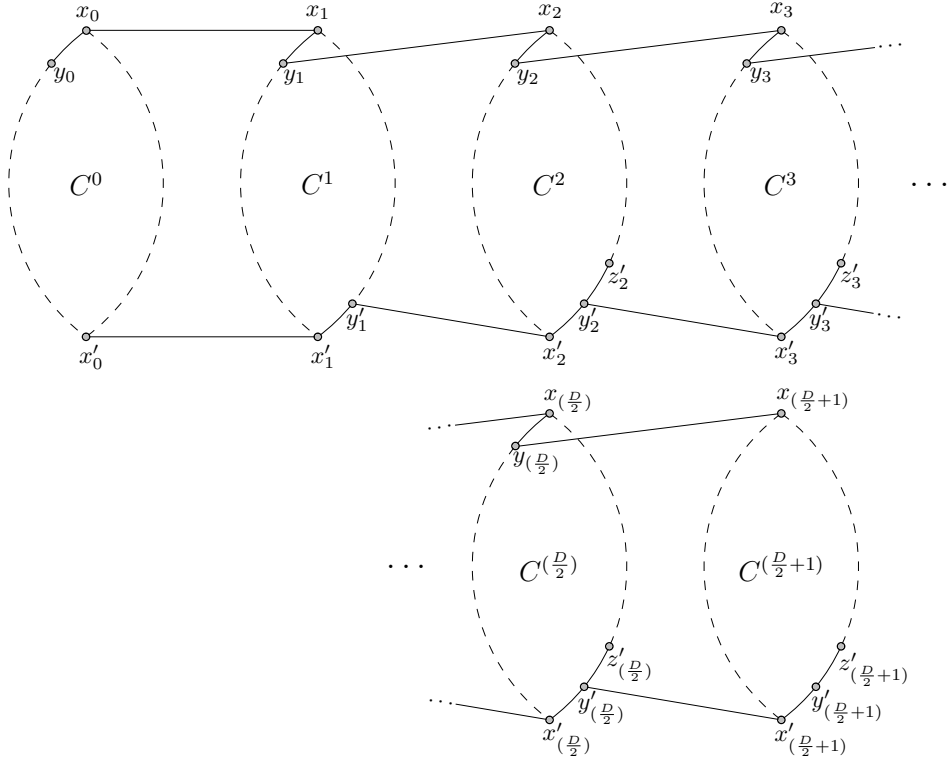}}
\caption{Auxiliary figure showing the sequence $C^0,C^1,\ldots,C^{D/2+1}$ of short cycles.}
\label{fig:Sequence}
\end{center}
\end{figure}

{\bf Claim 1.} $C^0\cap C=y_0x_0$.

{\bf Proof of Claim 1.} Let $y_0$, $z_0$, $y_2$, $y'_2$ and $z'_2$ be as in Fig.~\ref{fig:CubicEvenDiameter1} ($a$).

Consider a path $P=x'_2-y_0$. If $y'_1\in P$, then $P$ would go through a neighbour of $y'_1$ contained in $C^1$ and there would be a cycle in $\G$ of length at most $2D-4$. Therefore, we may assume $y'_2\in P$. If $\{x'_2,y'_2,z'_2\}\subset V(P\cap C2)$ then there would be a short cycle intersecting the cycle $C^2$ at a path of length $D-3$, a contradiction to Lemma \ref{lemm:D-3Int}. Similarly, we have that $z_0\not \in P$. Also, $P$ must be a $(D-1)$-path and $x_0\not\in P$; otherwise there would be a short cycle intersecting the cycle $C^2$ at a path of length $D-2$, a contradiction to Lemma \ref{lemm:(D-2)-Intersection}.

For $3\le i\le D/2+1$, let $x_i$ and $x'_i$ be the neighbours of $y_{i-1}$ and $y'_{i-1}$, respectively, not contained in  $C^{i-1}$, and let $C^i$ be (in virtue of the Saturating Lemma) the additional short cycle disjoint from $C^{i-1}$ which contains $x_i$ and $x'_i$. Since $P$ must go through $x'_i$, we denote by $y'_i$ the neighbour of $x'_i$  on $P\cap C^i$ and set $y_i=\rep^{C^i}(y'_i)$. We now show that, if $i\neq D/2+1$, $P\cap C^i=x'_iy'_i$. Assume the contrary; that is, $P\cap C^i=x'_iy'_iz'_i$ (since $\g(\G)=2D-2$, $|V(P\cap C^i)|\le3$). In such case, there would be a short cycle $y_0Pz'_iC^ix_iy_{i-1}x_{i-1}y_{i-2}x_{i-2}\ldots y_{1}x_{1}x_0y_0$ intersecting $C^i$ at a path of length $D-3$, contradicting Lemma \ref{lemm:D-3Int} (see Figures ~\ref{fig:CubicEvenDiameter1} ($b$) and ~\ref{fig:Sequence}). Consequently, $P\cap C^i=x'_iy'_i$ and $P$ must go through a neighbour of $y'_i$ not contained in $C^i$.

This way, for $3\le i\le D/2+1$ we have $d(y_0,y'_i)=d(y_0,y'_{i-1})-2=D-2(i-1)$, which means that $d(y_0,y'_{D/2+1})=0$. Since the cycle $C^{D/2+1}$ contains the vertices $y_0\in C^0\cap C\cap P$ and $x'_{D/2+1}\in P-C^0$, we have $C^{D/2+1}=C$, which implies that $C^0\cap C=y_0x_0$. \EndProof

As the selection of $C^0$ and $C$ was arbitrary, basically as a corollary of Claim 1 we have:

{\bf Claim 2.} Any two non-disjoint short cycles in $\G$ intersect at an edge.

\begin{figure}[!ht]
\begin{center}
\makebox[\textwidth][c]{\includegraphics[scale=.9]{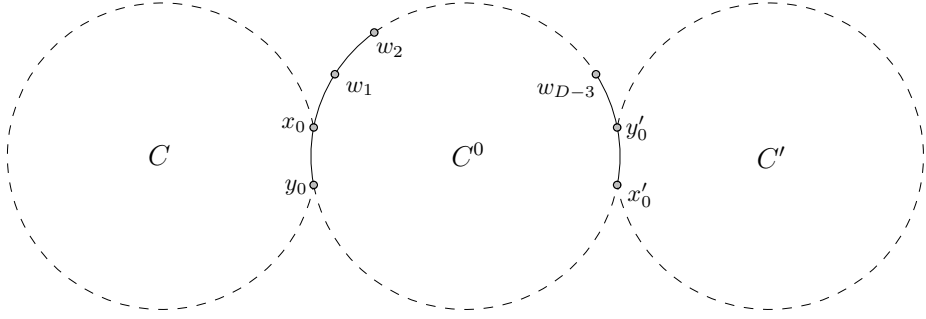}}
\caption{Auxiliary figure for Theorem \ref{theo:(3,D,-4)}.}
\label{TheoCubic}
\end{center}
\end{figure}

Finally, suppose that $C^0$ and $C$ intersect at $y_0x_0$, as stated by Claim 1. Let $y'_0$ be the repeat vertex of $y_0$ in $C^0$; then,  by Corollary \ref{cor:RepeatPath}, the repeat cycle $C'$ of $C$ intersects $C^0$ at $y'_0x'_0$ (the repeat path of $y_0x_0$ in $C^0$). Setting $Q=x_0C^0y'_0=x_0w_1\ldots w_{D-3}y'_0$, we have $Q$ is a path of length $D-2$ with saturated endvertices (see Fig.~\ref{TheoCubic}). Therefore, by Claim 2, there exists a sequence $F^1,\ldots,F^{D/2-2}$ of short cycles such that $F^i\cap C^0=w_{2i-1}w_{2i}$. However, since $D$ is even, the other short cycle containing $w_{D-3}$ would also contain one of the vertices in $\{w_{D-4},y'_0\}$, which contradicts the fact that $w_{D-4}$ and $y'_0$ are both saturated. This completes the proof of Theorem \ref{theo:(3,D,-4)}.\EndProof

Combining Theorems \ref{theo:OddDiameter} and \ref{theo:(3,D,-4)}, we have that the only bipartite $(3,D,-4)$-graphs with $D\ge 2$ are those depicted in Figures \ref{fig:Bipartite(3,3,-4)} and \ref{fig:Bipartite(3,4,-4)}, completing in this way the catalogue of such graphs.


\section{Conclusions}

The main results obtained in this paper are summarised below.

First we stated important structural properties of bipartite $(d,D,-4)$-graphs with $d\ge 3$ and $D\ge 3$. We found necessary conditions for the existence of bipartite $(d,D,-4)$-graphs with $d\ge 3$ and $D\ge 5$, which allowed us to prove the non-existence of such graphs for infinitely many pairs $(d,D)$; this included the case in which $D-1$ is an odd prime power, and the case in which $D-1\equiv 0 \pmod 3$. Afterwards,  we went on to proving that bipartite $(d,D,-4)$-graphs for $d\ge 3$ and odd $D\ge 5$ do not exist. We completed the catalogue of bipartite $(\D,D,-4)$-graphs with $\D\ge2$ and $5\le D\le 187$, which in turn completed the catalogue of bipartite $(\D,D,-\epsilon)$-graphs with $\D\ge2$, $5\le D\le 187$, $D\ne6$ and $0\le\epsilon\le4$.

\begin{description}
\item[Catalogue of bipartite $(\D,D,0)$-graphs with $\D\ge2$ and $5\le D\le 187$.] For $5\le D\le 187$ and $\D=2$ the only Moore bipartite graphs are the $2D$-cycles, whereas for $D=6$ and $\D\ge3$ they are incidence graphs of generalised polygons. For other values of $5\le D\le 187$ and $\D\ge 3$ there are no Moore bipartite graphs.

\item[Catalogue of bipartite $(\D,D,-2)$-graphs with $\D\ge2$ and $5\le D\le 187$.] The results of \cite{PV} combined with \cite{DJMP2,DJMP3} showed that there are no such graphs.

\item[Catalogue of bipartite $(\D,D,-4)$-graphs with $\D\ge2$ and $5\le D\le 187$.] The path of length 5 is the only such graph.
\end{description}

Another important result of the paper is the completion of the catalogue of bipartite $(3,D,-\epsilon)$-graphs with $D\ge 2$ and $0\le\epsilon\le4$.

\begin{description}
\item[Catalogue of bipartite $(3,D,0)$-graphs with $D\ge 2$.] The cubic Moore bipartite graphs are the complete bipartite graph $K_{3,3}$ for $D=2$, the unique incidence graph of the projective plane of order $2$ for $D=3$, the unique incidence graph of the generalised quadrangle of order $2$ for $D=4$, and the unique incidence graph of the generalised hexagon of order $2$ for $D=6$.

\item[Catalogue of  bipartite $(3,D,-2)$-graphs with $D\ge 2$.] There are only two non-isomorphic $(3,D,-2)$-graphs with $D\ge 2$; a unique bipartite $(3,2,-2)$-graph (the claw graph), and a unique $(3,3,-2)$-graph, which is depicted in Fig.~\ref{fig:Bipartite(d,3,-2)} ($a$).

\item[Catalogue of bipartite $(3,D,-4)$-graphs with $D\ge 2$.] There exist no bipartite $(3,2,-4)$-graphs. When the diameter is 3, there are four non-isomorphic bipartite $(3,3,-4)$-graphs; all of them are shown in Figure~\ref{fig:Bipartite(3,3,-4)}. For diameter 4, there is a unique bipartite $(3,4,-4)$-graph, which is depicted in Fig.~\ref{fig:Bipartite(3,4,-4)}.
    The results of this paper, combined with those of \cite{Jor93}, assert that there are no bipartite $(3,D,-4)$-graphs with $D\ge5$, outcome that gives an alternative proof of the optimality of the known bipartite $(3,5,-6)$-graph (see \cite{BD88}).
\end{description}

\subsection{Bipartite $(d,D,-4)$-graphs with $d\ge4$ and $D=3,4$}
\label{sub:Diameters34}

The main results in this paper did not include bipartite $(d,D,-4)$-graphs with $d\ge4$ and $D=3,4$. However, we believe that the structural properties of these graphs provided in Section \ref{sec:PreResults} could bear more conclusive results on such diameters.

For instance, by using Proposition \ref{prop:GirthBipartite(d,D,-4)}, Lemma \ref{lemm:SaturatingLemma} (Saturating Lemma), Lemma \ref{lemm:RepeatCycleLemma} (Repeat Cycle Lemma), Proposition \ref{prop:CyclePartition}, and Proposition \ref{prop:VertexPartition}, we were able to prove analytically the uniqueness of the two bipartite $(4,3,-4)$-graphs depicted in Fig.~\ref{fig:Bipartite(4,3,-4)}. We also think there should be no major difficulty to complete as well --in a very similar manner-- the catalogue of bipartite $(5,3,-4)$-graphs, which has so far as a unique element, the graph in Fig.~\ref{fig:Bipartite(5,3,-4)}.

Unfortunately, the final ideas used in the paper cannot be easily extended to cover bipartite $(d,D,-4)$-graphs with $d\ge4$ and $D=3,4$. With our current approach we cannot have Theorem \ref{theo:Theta} for $D=3,4$. In Theorem \ref{theo:Theta} the intersection paths $D^1\cap D^2$, $D^2\cap D^3$ and $D^3\cap D^4$ have length at most 2, and for us to apply the the Repeat Cycle Lemma we need the lengths of such paths to be  less than $D-2$. Indeed, the graph in Fig.~\ref{fig:Bipartite(4,3,-4)} ($a$) offers a good illustration of this. Even if we had Theorem \ref{theo:Theta}, something similar would occur with Lemma \ref{lemm:(D-2)-Intersection}; see the graphs in Fig.~\ref{fig:Bipartite(4,3,-4)} ($b$) and Fig.~\ref{fig:Bipartite(5,3,-4)}.

\subsection{Remarks on the upper bound for $\N^b(\D,D)$}


Our results improve the upper bound on $N^b_{\D,D}$ for many combinations of $\D$ and $D$. Recall that a bipartite $(\D,D,-5)$-graph $\G$ with $\D \ge 3$ and $D \ge 5$ must be regular (by Proposition \ref{prop:regularity1}) and thus cannot exist.

\begin{proposition}
For natural numbers $\D\ge 3$ and $D\ge8$ such that $D-1$ is an odd prime power, $\N^b(\D,D)\le \M^b(\D,D)-6$.
\end{proposition}

\begin{proposition}
For natural numbers $\D\ge 3$ and $D\ge7$ such that $d-1\equiv0,1\pmod p$, where $p$ is an odd prime factor of $D-1$, $\N^b(\D,D)\le \M^b(\D,D)-6$.
\end{proposition}

\begin{proposition}
For natural numbers  $\D\ge 3$ and $D\ge5$ such that $D\equiv1,3,4,5\pmod6$, $\N^b(\D,D)\le \M^b(\D,D)-6$.
\end{proposition}

\begin{proposition}
For natural numbers  $\D\ge 3$ and $5\le D\le187$ ($D\ne6$), $\N^b(\D,D)\le \M^b(\D,D)-6$.
\end{proposition}

\begin{proposition}
For any natural number $D\ge5$ ($D\ne6$), $\N^b(3,D)\le \M^b(3,D)-6$.
\end{proposition}

Finally, we feel that the next conjectures are valid.


\begin{conjecture}
There is no bipartite $(\D,D,-4)$-graph with $\D\ge 3$ and $D\ge 5$.
\end{conjecture}
\begin{conjecture}
For natural numbers $\D\ge 3$ and $D\ge 5$ such that $D\ne 6$, $\N^b(\D,D)\le \M^b(\D,D)-6$.
\end{conjecture}

\def\cprime{$'$}
\providecommand{\bysame}{\leavevmode\hbox to3em{\hrulefill}\thinspace}
\providecommand{\MR}{\relax\ifhmode\unskip\space\fi MR }
\providecommand{\MRhref}[2]{%
  \href{http://www.ams.org/mathscinet-getitem?mr=#1}{#2}
}
\providecommand{\href}[2]{#2}

\end{document}